\pdfoutput=1
\RequirePackage{ifpdf}
\ifpdf % We~are running pdfTeX in pdf mode
\documentclass[pdftex]{sigma}
\else
\documentclass{sigma}
\fi

\usepackage[all]{xy}

\newcommand{\Z}{\mathbb{Z}}
\newcommand{\N}{\mathbb{N}}

\newcommand{\id}{\mathbf{1}}

\numberwithin{equation}{section}

\newtheorem{Theorem}{Theorem}[section]
\newtheorem*{Theorem*}{Theorem}
\newtheorem{Corollary}[Theorem]{Corollary}
\newtheorem{Lemma}[Theorem]{Lemma}
\newtheorem{Proposition}[Theorem]{Proposition}
 { \theoremstyle{definition}
\newtheorem{Definition}[Theorem]{Definition}

\newtheorem{Example}[Theorem]{Example}
\newtheorem{Remark}[Theorem]{Remark} }

\begin{document}

\allowdisplaybreaks

\newcommand{\arXivNumber}{2401.15601}

\renewcommand{\PaperNumber}{080}

\FirstPageHeading

\ShortArticleName{On $F$-Polynomials for Generalized Quantum Cluster Algebras and Gupta's Formula}

\ArticleName{On $\boldsymbol{F}$-Polynomials for Generalized Quantum\\ Cluster Algebras and Gupta's Formula}

\Author{Changjian FU, Liangang PENG and Huihui YE}
\AuthorNameForHeading{C.~Fu, L.~Peng and H.~Ye}

\Address{Department of Mathematics, Sichuan University, Chengdu 610064, P.R.~China}
\Email{\href{mailto:changjianfu@scu.edu.cn}{changjianfu@scu.edu.cn}, \href{mailto:penglg@scu.edu.cn}{penglg@scu.edu.cn}, \href{mailto:yehuihuimath@outlook.com}{yehuihuimath@outlook.com}}

\ArticleDates{Received March 12, 2024, in final form August 25, 2024; Published online September 03, 2024}

\Abstract{We show the polynomial property of $F$-polynomials for generalized quantum cluster algebras and obtain the associated separation formulas under a mild condition. Along the way, we obtain Gupta's formulas of $F$-polynomials for generalized quantum cluster algebras. These formulas specialize to Gupta's formulas for quantum cluster algebras and cluster algebras respectively. Finally, a generalization of Gupta's formula has also been discussed in the setting of generalized cluster algebras.}


\Classification{13F60; 16S34; 05E16}

\section{Introduction}\label{S1}
Cluster algebras were invented by Fomin and Zelevinsky in \cite{fomin_zelevinsky_2002} with the aim to provide a combinatorial framework for the study of total positivity in algebraic groups and canonical bases of quantum groups. Since then, cluster algebras have been found deep connections with many other areas of mathematics and physics, such as discrete dynamical systems, non-commutative algebraic geometry, string theory and quiver representation theory etc., cf.\ \cite{Keller2012} and the references therein. A cluster algebra is a commutative algebra endowed with a distinguished set of generators called cluster variables. These generators are gathered into
overlapping sets of fixed finite cardinality, called clusters, which are defined recursively from an initial one via an operation called mutation. The first fundamental result in cluster algebras is the Laurent phenomenon~\cite{fomin_zelevinsky_2002}, which states that every cluster variable can be expressed as a Laurent polynomial in the initial ones. A basic problem in the structure theory of cluster algebras is to find an explicit expression of the Laurent polynomial of a cluster variable. Based on the Laurent phenomenon, Fomin and Zelevinsky \cite{fomin_zelevinsky_2007} further introduced $F$-polynomials and established the famous separation formulas, which give an expression of a cluster variable by its $g$-vector and $F$-polynomial. We remark that these separation formulas have played key roles not only in the structure theory of cluster algebras but also in the categorification of cluster algebras.\looseness=1

Quantum cluster algebras were introduced by Berenstein and Zelevinsky \cite{Berenstein_Zelevinsky2005}, which are $q$-deformations of cluster algebras of geometric type. It appears naturally in the study of algebraic varieties arising from Lie theory. Quantum cluster algebras share almost the same structure theory as the one of cluster algebras. Among others, Berenstein and Zelevinsky \cite{Berenstein_Zelevinsky2005} established the Laurent phenomenon for quantum cluster algebras. Tran \cite{Thao2011} further proved the existence of $F$-polynomials and established the separation formulas for quantum cluster algebras.\looseness=1

In their study of Teichm\"{u}ller space of Riemann surface with orbifold points, Chekhov and Shapiro \cite{Chekhov_Shapiro} discovered a new class of commutative algebras and formulated the so-called {\it generalized cluster algebras}. By definition, the notion is a generalization of cluster algebras. Chekhov and Shapiro \cite{Chekhov_Shapiro} found that the Laurent phenomenon is still true for generalized cluster algebras. Nakanishi \cite{Nakanishi14} proved that a generalized cluster algebra has almost the same structure theory of a cluster algebra. In particular, he introduced $F$-polynomials and established the separation formulas for generalized cluster algebras. It is worth mentioning that the structure of generalized cluster algebras also appear in many other branch of math, such as representation theory of quantum affine algebra \cite{Gleitz2015}, WKB analysis \cite{Iwaki-Nakanishi2016} and representation theory of finite-dimensional algebras \cite{LFM18,LF-Velasco2018}.

Inspired by quantum cluster algebras, it is natural to pursue a quantization of a generalized cluster algebra. However, it is not clear how to define a correct $q$-deformation of a generalized cluster algebra at this moment. Nevertheless, Nakanishi \cite{Nakanishi-qgca} considered the quantization of coefficients for generalized cluster algebras and discovered the quantum dilogarithms of higher degrees. Bai--Chen--Ding--Xu \cite{BCDX2018} introduced the notion of generalized quantum cluster algebra, which is a generalization of quantum cluster algebras and also a $q$-deformation of a %very
special class of generalized cluster algebras. Recently, the Laurent phenomenon for generalized quantum cluster algebras has been established in \cite{BCDX2023}.

The aim of this paper is to extend certain structure results of (quantum) cluster algebras to the setting of generalized (quantum) cluster algebras.
Under a mild condition on the mutation data $\mathbf{h}$, we prove the polynomial property %existence
of $F$-polynomials and establish the separation formulas for generalized quantum cluster algebras. The key ingredients of our approach are Fock--Goncharov decomposition of mutations and
Gupta type formulas.
As a byproduct, we obtain Gupta's formulas of $F$-polynomials for generalized quantum cluster algebras. These formulas degenerate to the Gupta's formulas of cluster algebras, which was discovered by Gupta~\cite{gupta2019} and recently reformulated and proved by Lin--Musiker--Nakanishi \cite{LMN2023}.
Since generalized quantum cluster algebras are not $q$-deformations of generalized cluster algebras with principal coefficients, the $F$-polynomials of a generalized quantum cluster algebra do not degenerate to $F$-polynomials of the associated generalized cluster algebra. Hence Gupta's formulas for generalized quantum cluster algebras do not degenerate to Gupta's formulas of generalized cluster algebras. Nevertheless, we show the strategy of \cite{LMN2023} can be extended to prove Gupta's formulas for generalized cluster algebras.

The paper is organized as follows. In Section \ref{Preliminaries}, we recollect basic results in generalized cluster algebras and generalized quantum cluster algebras. Section \ref{s:FG-decom-GQCA} is devoted to the study of Fock--Goncharov decomposition of mutations for generalized quantum cluster algebras. In Section \ref{s:F-poly_gen_quantum_cluster_alg}, we show the polynomial property of $F$-polynomials for generalized quantum cluster algebras and their associated Gupta's formulas. In Section \ref{s:gupta-formula-gca}, Gupta's formulas for generalized cluster algebras are discussed.

 \section{Preliminaries}\label{Preliminaries}
 In this section, we recall definitions and basic results for generalized cluster algebras \cite{Chekhov_Shapiro} and generalized quantum cluster algebras \cite{BCDX2018}.

 Throughout this section, we fix a positive integer $n$. Denote by $\mathbb{T}_n$ the $n$-regular tree whose edges are labeled by the numbers $1,\dots,n$ such that the $n$ edges emanating from each vertex carry different labels. We write
 \smash{$
 \xymatrix{t\ar@{-}[r]^k&t'}$}
  to indicate that the vertices $t$, $t'$ of $\mathbb{T}_n$ are linked by an edge labeled by $k$.
For an integer $b$, we use the notation $[b]_+:=\max(b,0)$. We also denote by $[1,n]$ the set $\{1,\dots,n\}$. Denote by $A^{\mathsf T}$ the transpose of a matrix $A$. A non-zero integer vector $\alpha\in \Z^n$ is {\it sign-coherent} if its entries are either all non-negative, or all non-positive.
For an integer vector $\alpha=(a_1,\dots,a_n)^{\mathsf T}\in\mathbb{Z}^n$, we denote by $[\alpha]_+:=([a_1]_+,\dots,[a_n]_+)^{\mathsf T}$. It is clear that $\alpha=[\alpha]_+-[-\alpha]_+$.

\subsection{Generalized cluster algebra}\label{S2}
 We follow \cite{Nakanishi14}.
Let $\mathbb{P}$ be a semifield, whose addition is denoted by $\oplus$. Denote by $\mathbb{ZP}$ the group ring of the multiplicative group $\mathbb{P}$ over $\mathbb{Z}$ and $\mathbb{QP}$ its fraction field. Let
$\mathcal{F}$ be the field of rational function in $n$ variables with coefficients in $\mathbb{QP}$.

\begin{Definition}\label{D4}
 A {\it labeled seed} with coefficients in $\mathbb{P}$ is a triplet $(\mathbf{x},\mathbf{y},B)$ such that
 \begin{itemize}\itemsep=0pt
 \item $B=(b_{ij})_{i,j=1}^{n}$ is a skew-symmetrizable integer matrix;
 \item $\mathbf{x}=(x_1,\dots,x_n)$ is an $n$-tuple of algebraic independent elements of $\mathcal{F}$ over $\mathbb{QP}$;
 \item $\mathbf{y}=(y_1,\dots,y_n)$ is an $n$-tuple of elements in $\mathbb{P}$.

 \end{itemize}
We say that $\mathbf{x}$ is a {\it cluster} and refer to $x_i$, $y_i$ and $B$ as the {\it cluster variables}, the {\it coefficients}
and the {\it exchange matrix}, respectively.
\end{Definition}
For a given labeled seed $(\mathbf{x},\mathbf{y},B)$ and $k\in [1,n]$, we set
\smash{$
\hat{y}_k:=y_k\prod_{j=1}^nx_j^{b_{jk}}
$}
and denote ${\hat{\mathbf{y}}=(\hat{y}_1,\dots, \hat{y}_n)}$. For an integer vector $\mathbf{a}=(a_1,\dots, a_n)^{\mathsf T}\in \Z^n$, we define
\[
\mathbf{x}^\mathbf{a}:=x_1^{a_1}\cdots x_n^{a_n},\qquad \mathbf{y}^\mathbf{a}:=y_1^{a_1}\cdots y_n^{a_n},\qquad \hat{\mathbf{y}}^\mathbf{a}:=\hat{y}_1^{a_1}\cdots \hat{y}_n^{a_n}.
\]

In order to introduce the mutation in generalized cluster algebra, we need the notion of mutation data.
\begin{Definition}\label{D5}
 A {\it mutation data} is a pair $(\mathbf{r},\mathbf{z})$, where
 \begin{itemize}\itemsep=0pt
 \item $\mathbf{r}=(r_1,\dots,r_n)$ is an $n$-tuple of positive integers;
 \item $\mathbf{z}=(z_{i,s})_{i=1,\dots,n;s=1,\dots,r_{i}-1}$ is a family of elements in $\mathbb{P}$ satisfying the reciprocity condition: $z_{i,s}=z_{i,r_i-s}$ for $1\leq s\leq r_{i}-1$.
 \end{itemize}
\end{Definition}
Throughout this subsection, we fix a mutation data $(\mathbf{r},\mathbf{z})$ and set $z_{i,0}=z_{i,r_i}=1$. Now we introduce the $(\mathbf{r},\mathbf{z})$-mutation in generalized cluster algebras.
\begin{Definition}\label{D6}
 For any seed ($\mathbf{x}$,$\mathbf{y}$,$B$) with coefficients in $\mathbb{P}$ and $k\in[1,n]$, the $(\mathbf{r},\mathbf{z})$-mutation of $(\mathbf{x},\mathbf{y},B)$ in direction $k$ is a new seed $\mu_k(\mathbf{x},\mathbf{y},B):=(\mathbf{x}',\mathbf{y}',B') $ with coefficients in $\mathbb{P}$ defined by the following rule:
\begin{align}
x_i'&=
\begin{cases}
 x_i & \text{if $i\neq k$}, \\
 \displaystyle x_k^{-1}\biggl(\prod\limits_{j=1}^{n}x_j^{[-\varepsilon b_{jk}]_+}\biggr)^{r_k}\frac{\sum\limits_{s=0}^{r_k}z_{k,s}\hat{y}_k^{\varepsilon s}}{\mathop{\bigoplus}\limits_{s=0}^{r_k}z_{k,s}y_k^{\varepsilon s}}& \text{if $i=k$},
\end{cases}\label{E4}
\\
 y_i'&=
\begin{cases}
 y_k^{-1} & \text{if $i= k$}, \\
\displaystyle y_i\bigl(y_k^{[\varepsilon b_{ki}]_+}\bigr)^{r_k}\biggl(\mathop{\bigoplus}\limits_{s=0}^{r_k}z_{k,s}y_k^{\varepsilon s}\biggr)^{-b_{ki}}& \text{if $i\neq k$},
\end{cases}\label{E5}
 \\
 b_{ij}'&=
 \begin{cases}
-b_{ij} &\text{if $i=k$ or $j=k$},\\
b_{ij}+r_k([-\varepsilon b_{ik}]_+b_{kj}+b_{ik}[\varepsilon b_{kj}]_+) &\text{else},
\end{cases}\label{E6}
\end{align}
where $\varepsilon\in\{\pm 1\}$.
\end{Definition}

\begin{Remark}\label{R1}\quad
\begin{itemize}\itemsep=0pt
\item[(1)] The mutation formulas \eqref{E4}, \eqref{E5} and \eqref{E6} are independent of the choice of $\varepsilon$ and $\mu_k$ is an involution.
 \item[(2)] If $\mathbf{r}=(1,\dots,1)$, then the mutation formulas \eqref{E4}, \eqref{E5} and \eqref{E6} reduce to the mutation formulas of cluster algebras.
 \item[(3)] The mutation of $\hat{\mathbf{y}}$ is similar as \eqref{E5} of $\mathbf{y}$:
 \begin{align*}
 \hat{y}_i'=
 \begin{cases}
 \hat{y}_k^{-1} & \text{if $i= k$}, \\
 \displaystyle \hat{y}_i\bigl(\hat{y}_k^{[\varepsilon b_{ki}]_+}\bigr)^{r_k}\left(\sum\limits_{s=0}^{r_k}z_{k,s}\hat{y}_k^{\varepsilon s}\right)^{-b_{ki}}& \text{if $i\neq k$}.
 \end{cases}
 \end{align*}
\end{itemize}
\end{Remark}

By assigning the labeled seed $(\mathbf{x},\mathbf{y}, B)$ to a root vertex $t_0\in \mathbb{T}_n$, we obtain an $(\mathbf{r},\mathbf{z})$-seed pattern $t\mapsto \Sigma_t$ of $(\mathbf{x},\mathbf{y}, B)$ in the same way as cluster algebras. In particular, for each vertex~${t\in \mathbb{T}_n}$, we have a labeled seed $\Sigma_t=(\mathbf{x}_t, \mathbf{y}_t, B_t)$ and if \smash{$\xymatrix{t\ar@{-}[r]^k &t'}$}, then $\Sigma_{t'}=\mu_k(\Sigma_t)$. We denote by $\mathbf{x}_t=(x_{1;t},\dots, x_{n;t})$, $\mathbf{y}_t=(y_{1;t},\dots,y_{n;t})$ and $B_t=(\mathbf{b}_{j;t})=(b_{ij;t})$.
\begin{Definition}
 The {\it generalized cluster algebra} $\mathcal{A}:=\mathcal{A}(t\mapsto \Sigma_t)$ associated to the $(\mathbf{r},\mathbf{z})$-seed pattern $t\mapsto \Sigma_t$ of $(\mathbf{x},\mathbf{y}, B)$ is the $\Z \mathbb{P}$-subalgebra of $\mathcal{F}$ generated by
 $\mathcal{X}:=\bigcup_{t\in \mathbb{T}_n}\mathbf{x}_t$.
\end{Definition}
The Laurent phenomenon still holds for generalized cluster algebras.
\begin{Proposition}[{\cite[Theorem 2.5]{Chekhov_Shapiro}}]\label{P6}
 Each cluster variable $x_{i;t}$ could be expressed as a Laurent polynomial of $\mathbf{x}$ with coefficients in $\mathbb{Z}\mathbb{P}$.
\end{Proposition}

\subsection{Generalized cluster algebra with principal coefficients}\label{ss:principal-gca}

From now on, let $\mathbf{y}=(y_1,\dots, y_n)$ and $\mathbf{z}=(z_{i,s})_{i=1,2,\dots,n;s=1,2,\dots,r_i-1}$ with $z_{i,s}=z_{i,r_{i}-s}$ be formal variables and $\mathbb{P}=\operatorname{Trop}(\mathbf{y},\mathbf{z})$ the tropical semifield of $\mathbf{y}$ and $\mathbf{z}$, which is the multiplicative abelian group freely
generated by $\mathbf{y}$ and $\mathbf{z}$ with tropical sum $\oplus$ defined by
\begin{equation*}
\biggl(\prod\limits_{i}y_i^{a_i}\prod\limits_{i,s}z_{i,s}^{a_{i,s}}\biggr)\oplus\biggl(\prod\limits_{i}y_i^{b_i}\prod\limits_{i,s}z_{i,s}^{b_{i,s}}\biggr)=\biggl(\prod\limits_{i}y_i^{\min\{a_i,b_i\}}\prod\limits_{i,s}z_{i,s}^{\min\{a_{i,s},b_{i,s}\}}\biggr),
\end{equation*}
where $a_i,a_{i,s},b_i,b_{i,s}\in \Z$. Let $(\mathbf{x},\mathbf{y}, B)$ be a labeled seed with coefficients in $\mathbb{P}$. Fix an $(\mathbf{r},\mathbf{z})$-seed pattern $t\mapsto \Sigma_t$ of $(\mathbf{x},\mathbf{y}, B)$ by assigning $(\mathbf{x},\mathbf{y}, B)$ to the vertex $t_0\in \mathbb{T}_n$. The associated generalized cluster algebra is called a {\it generalized cluster algebra with principal coefficients}. In this case, we denote it by $\mathcal{A}^\bullet$ to indicate the principal coefficients.

We assign two integer matrices $C_t=(\mathbf{c}_{1;t},\dots,\mathbf{c}_{n;t} )=(c_{ij;t})_{i,j=1}^n$ and $G_t=(\mathbf{g}_{1;t},\dots, \mathbf{g}_{n;t})=(g_{ij;t})_{i,j=1}^n$ to each vertex $t\in \mathbb{T}_n$ by the following recursion:
\begin{itemize}\itemsep=0pt
 \item $C_{t_0}=G_{t_0}=I_n$;
 \item if \smash{$\xymatrix{t\ar@{-}[r]^k&t'}\in \mathbb{T}_n$}, then
 \begin{align}
 c_{ij:t'}&=
 \begin{cases}
 -c_{ij;t} & \text{if $j=k$},\\
 c_{ij;t}+r_k(c_{ik;t}[\varepsilon b_{kj;t}]_++[-\varepsilon c_{ik;t}]_+b_{kj;t})& \text{if $ j\neq k$},
 \end{cases}\label{E9}\\
 \mathbf{g}_{i;t'}&=
 \begin{cases}
 \mathbf{g}_{i;t} & \text{if $i\neq k$},\\
 \displaystyle-\mathbf{g}_{k;t}+r_k\biggl(\sum\limits_{j=1}^{n}[-\varepsilon b_{jk;t}]_+\mathbf{g}_{j;t}-\sum\limits_{j=1}^{n}[-\varepsilon c_{jk;t}]_+\mathbf{b}_{j;t_0}\biggr) & \text{if $i=k$}.
 \end{cases}\label{E10}
 \end{align}
\end{itemize}
We remark that the recurrence formulas \eqref{E9} and \eqref{E10} are independent of the choice of the sign~${\varepsilon\in \{\pm 1\}}$.
We call $t\mapsto C_t$ and $t\mapsto G_t$ the~{\it $(\mathbf{r},\mathbf{z})$-$C$-pattern} and {\it $(\mathbf{r},\mathbf{z})$-$G$-pattern} of the~$(\mathbf{r},\mathbf{z})$-seed pattern of $(\mathbf{x},\mathbf{y}, B)$ respectively.
The column vectors of $C_t$ and $G_t$ are called {\it $c$-vectors} and~{\it $g$-vectors} of the $(\mathbf{r},\mathbf{z})$-seed pattern of $(\mathbf{x},\mathbf{y}, B)$, respectively. We remark that $t\mapsto C_t$ and~${t\mapsto G_t}$ only depend on $B, \mathbf{r}$ and $t_0\in \mathbb{T}_n$.

Denote by $R=\operatorname{diag}\{r_1,\dots,r_n\}$. Note that both $RB$ and $BR$ are skew-symmetrizable. Hence we may assign $(\mathbf{x},\mathbf{y}, RB)$ and $(\mathbf{x},\mathbf{y},BR)$ to the vertex $t_0$ to obtain (ordinary) seed patterns of~$(\mathbf{x},\mathbf{y}, RB)$ and $(\mathbf{x},\mathbf{y},BR)$ respectively.
It was proved by \cite[Proposition 3.9]{Nakanishi14} that the $C$-matrix $C_t$ of the $(\mathbf{r},\mathbf{z})$-seed pattern of $(\mathbf{x},\mathbf{y}, B)$ coincide with the $C$-matrix of the seed pattern of $(\mathbf{x},\mathbf{y}, RB)$ at vertex $t$. Alternatively, $R^{-1}C_tR$ is the $C$-matrix of the ordinary seed pattern of $(\mathbf{x},\mathbf{y}, BR)$ at vertex $t$. As a consequence, every $c$-vector of the $(\mathbf{r},\mathbf{z})$-seed pattern of $(\mathbf{x},\mathbf{y}, B)$ is sign-coherent. In this case, we also say $C_t$ is {\it column sign-coherent}. On the other hand, the $g$-vectors of the $(\mathbf{r},\mathbf{z})$-seed pattern of $(\mathbf{x},\mathbf{y}, B)$ at vertex $t$ coincide with the $g$-vectors of the ordinary seed pattern of $(\mathbf{x},\mathbf{y},BR)$ at vertex $t$. Hence the $G$-matrix $G_t$ is {\it row sign-coherent}, i.e., each row vector of $G_t$ is sign-coherent.
By the sign coherence of $c$-vectors provided in \cite[Corollary 5.5]{GHKK2018}, \eqref{E10} can be rewritten as
\begin{align}\label{E15}
\mathbf{g}_{i;t'}=\begin{cases}\mathbf{g}_{i;t}& \text{if $i\neq k$},\\
 \displaystyle-\mathbf{g}_{k;t}+r_k\biggl(\sum\limits_{j=1}^{n}[-\varepsilon_{k;t} b_{jk;t}]_+\mathbf{g}_{j;t}\biggr)& \text{if $i=k$},
 \end{cases}
\end{align}
where $\varepsilon_{k;t}$ is the common sign of components of the $c$-vector $\mathbf{c}_{k;t}$.

Similar to the ordinary seed pattern, we have the following tropical duality between $C$-matrices and $G$-matrices.
\begin{Proposition}[{\cite[Proposition 3.21]{Nakanishi14}}]\label{P4}
 Let $D_0$ be a diagonal matrix with positive integer diagonal entries such that $D_0RB$ is skew-symmetric. For each $t\in \mathbb{T}_n$, we have
\begin{align}\label{E11}
D_0^{-1}R^{-1}(G_t)^{\mathsf T}D_0RC_t=I_n.
\end{align}
\end{Proposition}
Let $D_0R=\operatorname{diag}\bigl\{d_1^{-1},\dots, d_n^{-1}\bigr\}$. We denote by $(-,-)_{D_0R}\colon \mathbb{Q}^n\times \mathbb{Q}^n\to \mathbb{Q}$ the inner product defined by $(\mathbf{u},\mathbf{v})_{D_0R}=\mathbf{u}^{\mathsf T}D_0R\mathbf{v}$, where $\mathbf{u},\mathbf{v}\in \mathbb{Q}^n$. With this notation, equation \eqref{E11} is equivalent to
\begin{align}\label{E12}
 (\mathbf{g}_{i;t},d_j\mathbf{c}_{j:t})_{D_0R}=\delta_{ij},\qquad \forall i,j\in[1,n].
\end{align}
Furthermore, by noticing that $D_0RB_t$ is skew-symmetric, we have
\begin{align*}
 (\mathbf{u},B_t\mathbf{v})_{D_0R}=-(B_t\mathbf{u},\mathbf{v})_{D_0R},\qquad \forall  \mathbf{u},\mathbf{v}\in \mathbb{Q}^n.
\end{align*}

\begin{Proposition}\label{P5}
 The following equality for the $(\mathbf{r},\mathbf{z})$-seed pattern $t\mapsto \Sigma_t$ holds:
 \begin{align}\label{E14}
 G_tB_t=B_{t_0}C_t.
 \end{align}
\end{Proposition}
\begin{proof}
Since $R^{-1}C_tR$ is the $C$-matrix of the seed pattern of $(\mathbf{x},\mathbf{y}, BR)$ and $G_t$ is the $G$-matrix of the seed pattern of $(\mathbf{x},\mathbf{y}, BR)$, we have
$
G_tB_tR=B_{t_0}R\bigl(R^{-1}C_tR\bigr)
$
by \cite[equation~(6.14)]{fomin_zelevinsky_2007}.
\end{proof}

For the generalized cluster algebra $\mathcal{A}^\bullet$ with principal coefficients, we have the strong Laurent phenomenon.
\begin{Proposition}[{\cite[Proposition 3.3]{Nakanishi14}}]
 Each cluster variable $x_{i;t}$ belongs to $\mathbb{Z}[\mathbf{x}^\pm,\mathbf{y},\mathbf{z}]$.
\end{Proposition}
\begin{Definition}
 The {\it $F$-polynomial} $F_{i;t}:=F_{i;t}[\mathbf{y},\mathbf{z}]$ of the cluster variable $x_{i;t}$ is defined as
 \[
 F_{i;t}[\mathbf{y},\mathbf{z}]:=x_{i;t}|_{x_1=\cdots=x_n=1}\in \mathbb{Z}[\mathbf{y},\mathbf{z}].
 \]
\end{Definition}

\begin{Proposition}[{\cite[Theorem 3.23]{Nakanishi14}}]\label{p:separation-formula-gca}
 For each $t\in \mathbb{T}_n$ and $i\in[1,n]$, the following formula holds:
 $
 x_{i;t}=\mathbf{x}^{\mathbf{g}_{i;t}}F_{i;t}[\hat{\mathbf{y}},\mathbf{z}]$.
\end{Proposition}
\begin{Remark}
 Similar to ordinary cluster algebras, the following statements are equivalent for generalized cluster algebras (cf.\ \cite[Proposition 3.19]{Nakanishi14}):
 \begin{itemize}\itemsep=0pt
 \item every {\it c-vectors} $\mathbf{c}_{i;t}$ is sign-coherent;
 \item every {\it F-polynomial} $F_{i;t}[\mathbf{y},\mathbf{z}]$ has a constant term 1;
 \item every {\it F-polynomial} $F_{i;t}[\mathbf{y},\mathbf{z}]$ has a unique monomial of maximal degree as a polynomial in~$\mathbf{y}$. Moreover its coefficients is 1 and it is divided by all the other occurring monomials as a polynomial in $\mathbf{y}$.
\end{itemize}
\end{Remark}

\subsection{Generalized quantum cluster algebras}\label{ss:gqca}
In this subsection, we introduce the definition and some properties of generalized quantum cluster algebras. We follow \cite{BCDX2018}.

Let $q$ be an indeterminate. Let $m\ge n$ be two positive integers, fix the mutation data $(R,\mathbf{h})$, where $R=\operatorname{diag}\{r_1,\dots,r_n\}$ is a diagonal $n\times n$ matrix whose diagonal coefficients are positive integers and $\mathbf{h}=(\mathbf{h}_1;\dots;\mathbf{h}_n)$ is defined as follows.
For $k\in[1,n]$,
\[
\mathbf{h}_{k}:=\bigl\{h_{k,0}\bigl(q^{\frac{1}{2}}\bigr),h_{k,1}\bigl(q^{\frac{1}{2}}\bigr),\dots,h_{k,r_{k}}\bigl(q^{\frac{1}{2}}\bigr)\bigr\},
\]
where $h_{k,i}\bigl(q^{\frac{1}{2}}\bigr)\in\mathbb{Z}\big[q^{\pm\frac{1}{2}}\big]$ satisfying $h_{k,i}\bigl(q^{\frac{1}{2}}\bigr)=h_{k,r_{k}-i}\bigl(q^{\frac{1}{2}}\bigr)$
and $h_{k,0}\bigl(q^{\frac{1}{2}}\bigr)=h_{k,r_{k}}\bigl(q^{\frac{1}{2}}\bigr)=1$.

A {\it compatible pair} $\bigl(\tilde{B}, \Lambda\bigr)$ consists of an integer $m\times n$-matrix $\tilde{B}$ and a skew-symmetric integer~${m\times m}$-matrix $\Lambda$ such that
\smash{$
\tilde{B}^{\mathsf T}\Lambda=[\begin{matrix} D &  0\end{matrix}]$},
where $D=\operatorname{diag}\bigl\{d_1^{-1},\dots, d_n^{-1}\bigr\}$ is a diagonal~${n\times n}$ matrix whose diagonal coefficients are positive integers. It is easy to see that the principal part~$B$ (i.e., the submatrix formed by the first $n$ rows) of $\tilde{B}$ is skew-symmetrizable and $D$ is a~skew-symmetrizer of $B$.

We define \smash{$E_{k,\varepsilon}^{\tilde{B}R}$} as the $m\times m$-matrix which differs from the identity matrix only in its $k$-th column whose coefficients are given by
 \begin{align*}
 \bigl(E_{k,\varepsilon}^{\tilde{B}R}\bigr)_{ik}=\begin{cases}
 -1& \text{if } i=k,\\
 [-\varepsilon b_{ik}r_k]_+& \text{if } i\neq k.
 \end{cases}
 \end{align*}
 Denote by \smash{$F_{k,\varepsilon}^{R\tilde{B}}$} the $n\times n$-matrix which differs from the identity matrix only in its $k$-th row whose coefficients are given by
 \begin{align*}
 \bigl(F^{R\tilde{B}}_{k,\varepsilon}\bigr)_{ki}=\begin{cases}
 -1& \text{if } i=k,\\
 [\varepsilon r_{k}b_{ki}]_+& \text{if } i\neq k.
 \end{cases}
 \end{align*}
Let $k\in [1,n]$.
The {\it mutation $\mu_k$ in direction $k$} transforms the compatible pair $\bigl(\tilde{B},\Lambda\bigr)$ into \smash{$\mu_{k}\bigl(\tilde{B}, \Lambda\bigr):=\bigl(\tilde{B}', \Lambda'\bigr)$}, where
\[
\tilde{B}'=E_{k,\varepsilon}^{\tilde{B}R}\tilde{B}F_{k,\varepsilon}^{R\tilde{B}},\qquad \Lambda'=\bigl(E_{k,\varepsilon}^{\tilde{B}R}\bigr)^{\mathsf T}\Lambda E_{k,\varepsilon}^{\tilde{B}R}.
\]
It is straightforward to check that the first equality is equivalent to \eqref{E6}. Moreover, $\bigl(\tilde{B}', \Lambda'\bigr)$ is a compatible pair and \smash{$\bigl(\tilde{B}'\bigr)^{\mathsf T}\Lambda'=[D\ 0]$}.

 Fix a skew-symmetric bilinear form $\lambda\colon\mathbb{Z}^{m}\times \mathbb{Z}^m\to \mathbb{Z}$. The {\it quantum torus $\mathcal{T}_\lambda$} associated~with~$\lambda$ is the $\mathbb{Z}\big[q^{\pm\frac{1}{2}}\big]$-algebra generated by the distinguished $\mathbb{Z}\big[q^{\pm\frac{1}{2}}\big]$-basis $\{\mathbf{x}(\alpha)\mid \alpha\in \mathbb{Z}^m\}$ with multiplication given by
\[
\mathbf{x}(\alpha)\mathbf{x}(\beta)=q^{\frac{1}{2}\lambda(\alpha,\beta)}\mathbf{x}(\alpha+\beta)
\]
for any $\alpha,\beta\in \mathbb{Z}^m$. The quantum torus $\mathcal{T}_\lambda$ is an Ore domain and we denote by $\mathcal{F}_q:=\mathcal{F}_\lambda$ its fraction skew field, which will be the ambient field to define the generalized quantum cluster algebras.
\begin{Definition}
 An {\it $(R,\mathbf{h})$-quantum seed} in $\mathcal{F}_q$ is a triple $\Sigma=\bigl(\mathbf{X},\tilde{B},\Lambda\bigr)$, where $\bigl(\tilde{B},\Lambda\bigr)$ is a~compatible pair and $\mathbf{X}=(X_1,\dots, X_m)$ is an $m$-tuple of elements of $\mathcal{F}_q$ such that
 \begin{itemize}\itemsep=0pt
 \item $X_1,\dots, X_{m}$ generated $\mathcal{F}_q$ over $\mathbb{Q}$;
 \item $X_iX_j=q^{\lambda_{ij}}X_jX_i$, where $\Lambda=(\lambda_{ij})$.
 \end{itemize}
 The (labeled) set $\mathbf{X}$ is called a {\it quantum cluster}, $X_1,\dots, X_n$ are {\it quantum cluster variables} and~${X_{n+1},\dots, X_m}$ are coefficients.
\end{Definition}
We define \[\mathbf{X}(\alpha):=q^{\frac{1}{2}\sum_{i<j}a_ia_j\lambda_{ji}}X_1^{a_1}\cdots X_{m}^{a_{m}}\] for any $\alpha=(a_1,\dots, a_{m})^{\mathsf T}\in \mathbb{Z}^{m}$. In particular, $X_i=\mathbf{X}(e_i)$, where $e_1,\dots, e_m$ is the standard $\Z$-basis of $\Z^m$.
It follows that
 \[
 \mathbf{X}(\alpha) \mathbf{X}(\beta)=q^{\frac{1}{2}\alpha^{\mathsf T}\Lambda \beta}\mathbf{X}({\alpha+\beta})
 \]
 for any $\alpha,\beta\in \Z^{m}$. The subalgebra of $\mathcal{F}_q$ generated by $\mathbf{X}(\alpha), \alpha\in \Z^{m}$, is a free $\Z\big[q^{\pm\frac{1}{2}}\big]$-module with basis $\{\mathbf{X}(\alpha) \mid \alpha\in \Z^{m}\}$, and hence identifies with the quantum torus $\mathcal{T}_\Lambda$ associated with the bilinear form $\Lambda\colon\Z^m\times \Z^m\to \Z, (\alpha,\beta)\mapsto \alpha^{\mathsf T}\Lambda \beta$, induced by $\Lambda$.

 For any $\beta\in \Z^n$, we also introduce the notation
$
 \hat{\mathbf{Y}}^{\beta}:=\mathbf{X}\bigl(\tilde{B}\beta\bigr)$.
 It is straightforward to check that
 \[
 \hat{\mathbf{Y}}^{\beta_1}\hat{\mathbf{Y}}^{\beta_2}=q^{\frac{1}{2}\beta_1^{\mathsf T}DB\beta_2}\hat{\mathbf{Y}}^{\beta_1+\beta_2},
 \]
 where $\beta_1,\beta_2\in \Z^n$.

 \begin{Definition}\label{mutation_of_generalzied_seed}
Let $\Sigma=\bigl(\mathbf{X},\tilde{B},\Lambda\bigr)$ be an $(R,\mathbf{h})$-quantum seed in $\mathcal{F}_q$. For any $k\in[1,n]$, the mutation $\mu_k$ in direction $k$ transforms the seed $\Sigma$ into a new triple \smash{$\mu_k(\Sigma):=\bigl(\mathbf{X}',\tilde{B}',\Lambda'\bigr)$}, where
\begin{itemize}\itemsep=0pt
\item[(1)] $\bigl(\tilde{B}', \Lambda'\bigr)=\mu_{k}\bigl(\tilde{B},\Lambda\bigr)$;
\item[(2)] $\mathbf{X}'=(X'_1,\dots, X'_m)$ is given by
\begin{align*}
 X'_i:=\mathbf{X}'(e_i)=\begin{cases}
 \mathbf{X}(e_i)&\text{if }i\neq k,\\
 \displaystyle\sum\limits_{s=0}^{r_k}h_{k,s}\bigl(q^{\frac{1}{2}}\bigr)\mathbf{X}(s[\varepsilon \mathbf{b}_k]_{+}+(r_{k}-s)[-\varepsilon \mathbf{b}_k]_+-e_i)&\text{if }i=k,
 \end{cases}
\end{align*}
where $\mathbf{b}_k$ is the $k$-th column vector of $\tilde{B}$ and $\varepsilon\in\{\pm 1\}$.
\end{itemize}
\end{Definition}
By \cite[Propositions 3.6 and 3.7]{BCDX2018}, the triple $\bigl(\mathbf{X}',\tilde{B}',\Lambda'\bigr)$ is also an $(R,\mathbf{h})$-quantum seed in $\mathcal{F}_q$ and $\mu_{k}$ is an involution.

 For a given $(R,\mathbf{h})$-quantum seed $\bigl(\mathbf{X}, \tilde{B}, \Lambda\bigr)$ in $\mathcal{F}_q$, we assign each vertex $t\in \mathbb{T}_n$ an $(R,\mathbf{h})$-quantum seed $\Sigma_t$ in $\mathcal{F}_q$ which can be obtained from \smash{$\bigl(\mathbf{X},\tilde{B},\Lambda\bigr)$} by iterated mutations such that if~$t$ and $t'$ are linked by an edge labeled $k$, then $\Sigma_{t'}=\mu_k(\Sigma_t)$. We call such an assignment $t\mapsto \Sigma_t$ an {\it $(R,\mathbf{h})$-quantum seed pattern}. It is clear that an $(R,\mathbf{h})$-quantum seed pattern is uniquely determined by the assignment of $\bigl(\mathbf{X},\tilde{B}, \Lambda\bigr)$ to an arbitrary vertex $t_0\in \mathbb{T}_n$. In this case, we refer to $t_0$ the {\it root vertex}, \smash{$\Sigma_{t_0}=\bigl(\mathbf{X}, \tilde{B},\Lambda\bigr)$} the {\it initial $(R,\mathbf{h})$-quantum seed} and $X_1,\dots,X_n$ the {\it initial quantum cluster variables}.
 In the following, when we fix an $(R,\mathbf{h})$-quantum seed pattern, we always denote by $\Sigma_t=\bigl(\mathbf{X}_t,\tilde{B}_t,\Lambda_t\bigr)$ and
\[\mathbf{X}_t=(X_{1;t},\dots,X_{m;t}),\qquad \tilde{B}_t=(b_{ij;t}),\qquad \Lambda_t=(\lambda_{ij;t}).
\]
For the initial quantum seed $\Sigma_{t_0}$, we denote
\[
\mathbf{X}_{t_0}:=\mathbf{X}=(X_1,\dots, X_{m}),\qquad \tilde{B}_{t_0}:=\tilde{B}=(b_{ij}), \qquad\Lambda_{t_0}:=\Lambda=(\lambda_{ij}).
\]
For each vertex $t$, we refer to $\mathbf{X}_t$ a {\it quantum cluster}, $X_{i;t}$ $(1\leq i\leq n)$ {\it quantum cluster variables} and $X_{n+i;t}$ $(1\leq i\leq m-n)$ {\it coefficients}.
It is clear that for every $t\in \mathbb{T}_n$, we have $X_{n+i;t}=X_{n+i}$ for $1\leq i\leq m-n$.
\begin{Definition}
 Given an $(R,\mathbf{h})$-quantum seed pattern $t\mapsto \Sigma_t$, the $(R,\mathbf{h})$-quantum cluster algebra $\mathcal{A}_q(\Sigma_{t_0})$ with initial seed $\Sigma_{t_0}=\bigl(\mathbf{X},\tilde{B},\Lambda\bigr)$ is the \smash{$\mathbb{Z}\big[q^{\pm\frac{1}{2}}\big]\big[X_{n+1}^{\pm 1},\dots,X_m^{\pm 1}\big]$} subalgebra of~$\mathcal{F}_q$ generated by all the quantum cluster variables
 \[\mathcal{X}_q:=\{X_{i;t}\mid 1\le i\le n, t\in \mathbb{T}_n\}.\]
\end{Definition}
\begin{Remark}
The algebra $\mathcal{A}_q(\Sigma_{t_0})$ is also called the {\it generalized quantum cluster algebra} associated with the $(R,\mathbf{h})$-quantum seed $\bigl(\mathbf{X},\tilde{B},\Lambda\bigr)$. In particular, if $R=\operatorname{diag}\{1,\dots, 1\}$, then~$\mathcal{A}_q(\Sigma_{t_0})$ is a quantum cluster algebra in the sense of Berenstein and Zelevinsky \cite{Berenstein_Zelevinsky2005}. On the other hand, it can be viewed as a $q$-deformation for a special generalized cluster algebra in the sense of Chekhov and Shapiro \cite{Chekhov_Shapiro}.
\end{Remark}

Fix an $(R,\mathbf{h})$-quantum seed pattern $t\mapsto \Sigma_t$ with initial seed $\Sigma_{t_0}=\bigl(\mathbf{X},\tilde{B},\Lambda\bigr)$. Denote by~${\mathbf{r}=(r_1,\dots, r_n)}$. Recall that we have the $(\mathbf{r},\mathbf{z})$-$C$-pattern $t\mapsto C_t$ and $(\mathbf{r},\mathbf{z})$-$G$-pattern $t\mapsto G_t$ associated with $B$, $\mathbf{r}$ and $t_0\in \mathbb{T}_n$ as in Section \ref{ss:principal-gca}.

We introduce the $(R,\mathbf{h})$-$G$-pattern $t\mapsto \tilde{G}_t$ for the $(R,\mathbf{h})$-quantum seed pattern $t\mapsto \Sigma_t$ as follows.
 For each vertex $t\in \mathbb{T}_n$, we assign an $m\times m$-integer matrix $\tilde{G}_t=(\tilde{\mathbf{g}}_{1;t},\dots, \tilde{\mathbf{g}}_{m;t})$ to $t$ by the following recursion:
\begin{itemize}\itemsep=0pt
 \item[(1)] $\tilde{G}_{t_0}=I_m$;
 \item[(2)] if \smash{$\xymatrix{t\ar@{-}[r]^k&t'}\in \mathbb{T}_n$}, then
 \begin{align}\label{eq:g-formula}
 \mathbf{\tilde{g}}_{i;t'}=
 \begin{cases}
 \mathbf{\tilde{g}}_{i;t} & \text{if } i\neq k,\\
 \displaystyle -\mathbf{\tilde{g}}_{k;t}+r_{k}\biggl(\sum\limits_{j=1}^{m}[-b_{jk;t}]_+\mathbf{\tilde{g}}_{j;t}-\sum\limits_{j=1}^{n}[- c_{jk;t}]_+\mathbf{b}_{j;t_0}\biggr) & \text{if }i=k.
 \end{cases}
 \end{align}

\end{itemize}
We call $\tilde{G}_t$ the {\it $G$-matrix} of the $(R,\mathbf{h})$-quantum seed pattern $t\mapsto \Sigma_t$ and the column vector $\tilde{\mathbf{g}}_{i;t}$ the {\it $g$-vector} of the quantum cluster variable $X_{i;t}$.
Clearly, we have \smash{$\tilde{G}_t=\left[\begin{smallmatrix}
 G_t&0\\ \star& I_{m-n}
\end{smallmatrix}\right]$}. It follows that
\begin{align}\label{D_form_duality}
 \mathbf{c}^{\mathsf T}_{i;t}[D\quad 0]\mathbf{\tilde{g}}_{j;t}=(\mathbf{c}_{i;t},\mathbf{g}_{j;t})_D=d_{i}^{-1}\delta_{ij}
\end{align}
for $i,j\in[1,n]$.

Similar to \cite[equation (6.14)]{fomin_zelevinsky_2007}, we have the following.
\begin{Proposition}
 For each vertex $t\in \mathbb{T}_n$, the following equation holds:
 \begin{align}\label{duality_quantum_G_C_matrix1}
 \tilde{G}_{t}\tilde{B}_{t}=\tilde{B}_{t_0}C_{t}.
 \end{align}
\end{Proposition}
\begin{proof}

 We prove the equality by induction on the distance between $t_0$ and $t$.
 It is obvious for~${t=t_0}$. Let \smash{$\xymatrix{t\ar@{-}[r]^{k}&t'}$} be an edge in $\mathbb{T}_n$ and suppose that \eqref{duality_quantum_G_C_matrix1} holds for the vertex $t$. We check it for $t'$. We first consider the $k$-th column,
 \begin{align*}
 \begin{aligned}
 \sum\limits^{m}_{i=1}b_{ik;t'}\mathbf{\tilde{g}}_{i;t'}
 &=\sum\limits_{i\neq k}(-b_{ik;t})\mathbf{\tilde{g}}_{i;t}=-\sum\limits^{m}_{i=1}b_{ik;t}\mathbf{\tilde{g}}_{i;t}=-\sum\limits^{n}_{i=1}c_{ik;t}\mathbf{b}_{i;t_0}\qquad
 \text{(by induction)}\\
 &=\sum\limits^{n}_{i=1}c_{ik;t'}\mathbf{b}_{i;t_0}.
 \end{aligned}
 \end{align*}
 Now let $j\neq k$, we have
 \begin{align*}
 \sum\limits^{m}_{i=1}b_{ij;t'}\mathbf{\tilde{g}}_{i;t'}
 ={}&\sum\limits_{i\neq k}b_{ij;t'}\mathbf{\tilde{g}}_{i;t'}+b_{kj;t'}\mathbf{\tilde{g}}_{k;t'}\\
 ={}&\sum\limits_{i\neq k}(b_{ij;t}+r_{k}(b_{ik;t}[b_{kj;t}]_{+}+[-b_{ik;t}]_{+}b_{kj;t}))\mathbf{\tilde{g}}_{i;t}\\
 & +(-b_{kj;t})\left(-\mathbf{\tilde{g}}_{k;t}+r_{k}\left(\sum\limits^{m}_{l=1}[-b_{lk;t}]_{+}\mathbf{\tilde{g}}_{l;t}-\sum\limits^{n}_{l=1}[-c_{lk;t}]_{+}\mathbf{b}_{l;t_0}\right)\right)\\
 ={}&\sum\limits^{m}_{i=1}b_{ij;t}\mathbf{\tilde{g}}_{i;t}+r_{k}[b_{kj;t}]_{+}\sum\limits^{m}_{i=1}b_{ik;t}\mathbf{\tilde{g}}_{i;t}+r_{k}b_{kj;t}\sum\limits^{n}_{i=1}[-c_{ik;t}]_{+}\mathbf{b}_{i;t_0}\\
 ={}&\sum\limits_{i=1}^{n}c_{ij;t}\mathbf{b}_{i;t_0}+r_{k}[b_{kj;t}]_{+}\sum\limits^{n}_{i=1}c_{ik;t}\mathbf{b}_{i;t_0}+r_{k}b_{kj;t}\sum\limits^{n}_{i=1}[-c_{ik;t}]_{+}\mathbf{b}_{i;t_0}\\
 ={}&\sum\limits_{i=1}^{n}(c_{ij;t}+r_{k}(c_{ik;t}[b_{kj;t}]_{+}+[-c_{ik;t}]_{+}b_{kj;t}))\mathbf{b}_{i;t_0}\\
 ={}&\sum\limits_{i=1}^{n}c_{ij;t'}\mathbf{b}_{i;t_0}.
 \end{align*}
 Thus \eqref{duality_quantum_G_C_matrix1} holds for the vertex $t'$. This completes the proof.
\end{proof}

By \eqref{duality_quantum_G_C_matrix1}, formula \eqref{eq:g-formula} is equivalent to the following:
 \begin{align}\label{eq:g-formula-new}
 \mathbf{\tilde{g}}_{i;t'}=
 \begin{cases}
 \mathbf{\tilde{g}}_{i;t} & \text{if } i\neq k,\\
 \displaystyle -\mathbf{\tilde{g}}_{k;t}+r_{k}\biggl(\sum\limits_{j=1}^{m}[-\varepsilon b_{jk;t}]_+\mathbf{\tilde{g}}_{j;t}-\sum\limits_{j=1}^{n}[-\varepsilon c_{jk;t}]_+\mathbf{b}_{j;t_0}\biggr) & \text{if }i=k,
 \end{cases}
 \end{align}
 where $\varepsilon\in\{\pm 1\}$.
By the sign-coherence of $c$-vectors, we have
\begin{align}\label{eq:g-matrix}
 \tilde{G}_{t'}=\tilde{G}_{t}E_{k,\varepsilon_{k;t}}^{\tilde{B}_{t}R},
\end{align}
where $\varepsilon_{k;t}$ is common sign of components of the $c$-vector $\mathbf{c}_{k;t}$.

If $m=2n$ and $\tilde{B}=\left[\begin{smallmatrix}
 B\\ I_n
 \end{smallmatrix}\right]$, then $(R,\mathbf{h})$-quantum cluster algebra $\mathcal{A}_q$ is called a {\it generalized quantum cluster algebra with principal coefficients}. In this case, the $g$-vector $\tilde{\mathbf{g}}_{i;t}$ is closely related to the $g$-vector $\mathbf{g}_{i;t}$.
\begin{Lemma}\label{l:g-vector-comparison}
 Let $m=2n$ and $\tilde{B}=\left[\begin{smallmatrix}
 B\\ I_n
 \end{smallmatrix}\right]$. Then for each vertex $t\in \mathbb{T}_n$ and $1\leq i\leq n$, we have~${\tilde{\mathbf{g}}_{i;t}=\left[\begin{smallmatrix}
 \mathbf{g}_{i;t}\\ 0
 \end{smallmatrix}\right]}$.
\end{Lemma}
\begin{proof}
 Note that $b_{jk;t}=c_{(j-n)k;t}$ whenever $n<j\leq 2n$. By the sign-coherence of $c$-vectors, we may choose a sign $\varepsilon$ such that $[-\varepsilon c_{(n-j)k;t}]_+=0$. Hence the recurrence formula \eqref{eq:g-formula-new} can be rewritten as
 \begin{align*}
 \tilde{\mathbf{g}}_{i;t'}= \begin{cases}
 \mathbf{\tilde{g}}_{i;t} & \text{if } i\neq k,\\
 \displaystyle-\mathbf{\tilde{g}}_{k;t}+r_{k}\sum\limits_{j=1}^{n}[-\varepsilon b_{jk;t}]_+\mathbf{\tilde{g}}_{j;t} & \text{if }i=k.
 \end{cases}
 \end{align*}
 Now the result can be deduced by induction on the distance between $t$ and $t_0$.
\end{proof}

\section[F-polynomials for generalized quantum cluster algebras]{$\boldsymbol{F}$-polynomials for generalized quantum cluster algebras}

Throughout this section, fix an $(R,\mathbf{h})$-quantum seed pattern $t\mapsto \Sigma_t$ with initial seed $\Sigma_{t_0}=\bigl(\mathbf{X},\tilde{B},\Lambda\bigr)$ and keep the notation in Section \ref{ss:gqca}. Recall that $e_1,\dots, e_m$ is the standard $\mathbb{Z}$-basis of $\mathbb{Z}^m$. We denote by $f_1,\dots, f_n$ the standard $\mathbb{Z}$-basis of $\mathbb{Z}^n$.
For an integer vector $\beta=(b_1,\dots, b_m)^{\mathsf T}\in \Z^m$, we denote by $\overline{\beta}=(b_1,\dots, b_n)^{\mathsf T}\in \Z^n$ the truncation of $\beta$.
\subsection{Fock--Goncharov decomposition}\label{s:FG-decom-GQCA}
In this section, we introduce Fock--Goncharov decomposition of mutations for generalized quantum cluster algebras, which generalizes the corresponding construction of quantum cluster algebras in \cite{fock2009,Keller2012}.

For $a,b\in\mathbb{Z}$ and $k\in[1,n]$, let
\begin{align*}
 \left(\sum\limits^{r_{k}}_{s=0}h_{k,s}\bigl(q^{\frac{1}{2}}\bigr)\bigl(q^{\frac{b}{2}}z\bigr)^{s}\right)^{\{a\}}:=
 \begin{cases}\displaystyle\prod\limits_{i=1}^{a}
 \left(\sum\limits^{r_{k}}_{s=0}h_{k,s}\bigl(q^{\frac{1}{2}}\bigr)\bigl(q^{\frac{b(2i-1)}{2}}z\bigr)^{s}\right)&\text{if } a>0,\\
 1&\text{if }a=0,\\
 \displaystyle \prod\limits_{i=a}^{-1}\left(\sum\limits^{r_{k}}_{s=0}h_{k,s}\bigl(q^{\frac{1}{2}}\bigr)\bigl(q^{\frac{b(2i+1)}{2}}z\bigr)^{s}\right)^{-1}&\text{if } a<0.
 \end{cases}
\end{align*}
It is easy to verify that
\begin{gather}
 \left(\sum\limits^{r_{k}}_{s=0}h_{k,s}\bigl(q^{\frac{1}{2}}\bigr)\bigl(q^{\frac{b}{2}}z\bigr)^{s}\right)^{\{a+a'\}}\nonumber\\
 \qquad=\left(\sum\limits^{r_{k}}_{s=0}h_{k,s}\bigl(q^{\frac{1}{2}}\bigr)\bigl(q^{\frac{b}{2}}x\bigr)^{s}\right)^{\{a\}}\bigg|_{x=q^{a'b}z}
 \left(\sum\limits^{r_{k}}_{s=0}h_{k,s}\bigl(q^{\frac{1}{2}}\bigr)\bigl(q^{\frac{b}{2}}z\bigr)^{s}\right)^{\{a'\}}.\label{eq:FG-decom}
\end{gather}

Let $\mathcal{T}_{t}$ be the quantum torus associated with $\Lambda_{t}$, i.e., the $\Z\big[q^{\pm\frac{1}{2}}\big]$-subalgebra of $\mathcal{F}_q$ generated by $\{\mathbf{X}_t(\alpha)\mid \alpha\in \Z^m\}$. It is an Ore domain. Denote by $\mathcal{F}_t$ the fraction skew field of $\mathcal{T}_t$.
Let~\smash{$\xymatrix{t\ar@{-}[r]^k&t'}$} be an edge in $\mathbb{T}_n$. The mutation $\mu_k$ in direction $k$ yields a unique \smash{$\Z\big[q^{\pm\frac{1}{2}}\big]$}-algebra isomorphism $\mu_{k;t}\colon\mathcal{F}_{t'}\to \mathcal{F}_t$ such that
\[
 \mu_{k;t}(\mathbf{X}_{t'}(e_i))=\begin{cases}
 \mathbf{X}_t(e_i) &\text{if $i\neq k$,}\\
 \displaystyle \sum\limits_{s=0}^{r_k}h_{k,s}\bigl(q^{\frac{1}{2}}\bigr)\mathbf{X}_t(s[\varepsilon \mathbf{b}_k]_{+}+(r_{k}-s)[-\varepsilon \mathbf{b}_k]_+-e_i)&\text{if $i=k$.}
 \end{cases}
\]

Recall that for $\alpha\in \Z^n$, we have $\hat{\mathbf{Y}}_t^\alpha=\mathbf{X}_t\bigl(\tilde{B}_t\alpha\bigr)$. For $k\in [1,n]$, we also denote \smash{$\hat{\mathbf{Y}}_{k;t}:=\mathbf{\hat{Y}}_t^{f_k}$}.
\begin{Lemma}
 For $k\in[1,n],t\in\mathbb{T}_{n}$ and $\alpha\in\mathbb{Z}^{n}$, there is a unique $\Z\big[q^{\pm\frac{1}{2}}\big]$-algebra homomorphism~\smash{$\psi_{k;t}\bigl(\mathbf{\hat{Y}}_{t}^{\alpha}\bigr)\colon\mathcal{F}_t\to \mathcal{F}_t$} such that
 \[
\psi_{k;t}\bigl(\mathbf{\hat{Y}}_{t}^{\alpha}\bigr)(\mathbf{X}_t(\beta))=\mathbf{X}_{t}(\beta)\left(\sum\limits_{s=0}^{r_k} h_{k,s}\bigl(q^{\frac{1}{2}}\bigr)\bigl(q^{\frac{1}{2d_{k}}}\mathbf{\hat{Y}}_{t}^{\alpha}\bigr)^{s}\right)^{-\{d_{k}(\bar{\beta},\alpha)_D\}}, \qquad \forall \beta\in \mathbb{Z}^m.
 \]
 \end{Lemma}

 \begin{proof}
It suffices to show that
 \[\psi_{k;t}\bigl(\mathbf{\hat{Y}}_{t}^{\alpha}\bigr)(\mathbf{X}_{t}(\beta_1)\mathbf{X}_{t}(\beta_2))=\psi_{k;t}\bigl(\mathbf{\hat{Y}}_{t}^{\alpha}\bigr)(\mathbf{X}_{t}(\beta_1))\psi_{k;t}\bigl(\mathbf{\hat{Y}}_{t}^{\alpha}\bigr)(\mathbf{X}_{t}(\beta_2))\]
for $\beta_{1},\beta_2\in\mathbb{Z}^{m}$. We have
\begin{align*}
 \psi_{k;t}\bigl(\mathbf{\hat{Y}}_{t}^{\alpha}\bigr)(\mathbf{X}_{t}(\beta_1)\mathbf{X}_{t}(\beta_2))
 ={}&\psi_{k;t}\bigl(\mathbf{\hat{Y}}_{t}^{\alpha}\bigr)\bigl(q^{\frac{1}{2}\beta_{1}^{t}\Lambda_{t}\beta_{2}}\mathbf{X}_{t}(\beta_1+\beta_2)\bigr)\\
 ={}&q^{\frac{1}{2}\beta_{1}^{t}\Lambda_{t}\beta_{2}}\mathbf{X}_{t}(\beta_{1}+\beta_{2})\left(\sum\limits_{s=0}^{r_k} h_{k,s}\bigl(q^{\frac{1}{2}}\bigr)\bigl(q^{\frac{1}{2d_{k}}}\mathbf{\hat{Y}}_{t}^{\alpha}\bigr)^{s}\right)^{-\{d_{k}(\bar{\beta_{1}}+\bar{\beta_{2}},\alpha)_D\}}\\
 ={}&\mathbf{X}_{t}(\beta_{1})\mathbf{X}_{t}(\beta_{2})\left(\sum\limits_{s=0}^{r_k} h_{k,s}\bigl(q^{\frac{1}{2}}\bigr)\bigl(q^{\frac{1}{2d_{k}}}\mathbf{\hat{Y}}_{t}^{\alpha}\bigr)^{s}\right)^{-\{d_{k}(\bar{\beta_{1}},\alpha)+d_k(\bar{\beta_{2}},\alpha)_D\}}\\
={}&\mathbf{X}_{t}(\beta_{1})\mathbf{X}_{t}(\beta_{2})\left(\sum\limits_{s=0}^{r_k} h_{k,s}\bigl(q^{\frac{1}{2}}\bigr)(q^{\frac{1}{2d_{k}}+ (\bar{\beta_{2}},\alpha)_D}\mathbf{\hat{Y}}_{t}^{\alpha})^{s}\right)^{-\{d_{k}(\bar{\beta_{1}},\alpha)_D\}}\\
&\times\left(\sum\limits_{s=0}^{r_k} h_{k,s}\bigl(q^{\frac{1}{2}}\bigr)\bigl(q^{\frac{1}{2d_{k}}}\mathbf{\hat{Y}}_{t}^{\alpha}\bigr)^{s}\right)^{-\{d_{k}(\bar{\beta_{2}},\alpha)_D\}}\\
={}&\mathbf{X}_{t}(\beta_{1})\left(\sum\limits_{s=0}^{r_k} h_{k,s}\bigl(q^{\frac{1}{2}}\bigr)\bigl(q^{\frac{1}{2d_{k}}}\mathbf{\hat{Y}}_{t}^{\alpha}\bigr)^{s}\right)^{-\{d_{k}(\bar{\beta_{1}},\alpha)_D\}}\\
&\times\mathbf{X}_{t}(\beta_{2})\left(\sum\limits_{s=0}^{r_k} h_{k,s}\bigl(q^{\frac{1}{2}}\bigr)\bigl(q^{\frac{1}{2d_{k}}}\mathbf{\hat{Y}}_{t}^{\alpha}\bigr)^{s}\right)^{-\{d_{k}(\bar{\beta_{2}},\alpha)_D\}}\\
={}&\psi_{k;t}\bigl(\mathbf{\hat{Y}}_{t}^{\alpha}\bigr)(\mathbf{X}_{t}(\beta_1))\psi_{k;t}\bigl(\mathbf{\hat{Y}}_{t}^{\alpha}\bigr)(\mathbf{X}_{t}(\beta_2)),
\end{align*}
where the fourth equality follows from equation \eqref{eq:FG-decom}.
\end{proof}

\begin{Lemma}
 For $\alpha\in \Z^n$, $\psi_{k;t}\bigl(\mathbf{\hat{Y}}_{t}^{\alpha}\bigr)$ is an isomorphism. Moreover, its inverse is given by
 \begin{align*}
 \psi_{k;t}\bigl(\mathbf{\hat{Y}}_{t}^{\alpha}\bigr)^{-1}\colon\ \mathcal{F}_t\to \mathcal{F}_t,\qquad
 \mathbf{X}_{t}(\beta)\mapsto\mathbf{X}_{t}(\beta)\left(\sum\limits_{s=0}^{r_k} h_{k,s}\bigl(q^{\frac{1}{2}}\bigr)\bigl(q^{\frac{1}{2d_{k}}}\mathbf{\hat{Y}}_{t}^{\alpha}\bigr)^{s}\right)^{\{d_{k}(\bar{\beta},\alpha)_D\}}.
\end{align*}
\end{Lemma}
\begin{proof}
 It is straightforward to show that there is a unique $\Z\big[q^{\pm\frac{1}{2}}\big]$-algebra homomorphism $\Phi\colon \mathcal{F}_t\to \mathcal{F}_t$ which maps $\mathbf{X}_t(\beta)$ to
 \[
 \mathbf{X}_{t}(\beta)\left(\sum\limits_{s=0}^{r_k} h_{k,s}\bigl(q^{\frac{1}{2}}\bigr)\bigl(q^{\frac{1}{2d_{k}}}\mathbf{\hat{Y}}_{t}^{\alpha}\bigr)^{s}\right)^{\{d_{k}(\bar{\beta},\alpha)_D\}}
 \]
  for any $\beta\in \Z^m$ by \eqref{eq:FG-decom}.
 Furthermore, one can show that
 \[
 \Phi\circ \psi_{k;t}\bigl(\mathbf{\hat{Y}}_{t}^{\alpha}\bigr)(\mathbf{X}_t(\beta))=\mathbf{X}_t(\beta)=\psi_{k;t}\bigl(\mathbf{\hat{Y}}_{t}^{\alpha}\bigr)\circ \Phi(\mathbf{X}_t(\beta))
 \]
 by noticing that $(B_t\alpha, \alpha)_D=0$, which implies the result.
\end{proof}

\begin{Proposition}\label{quantum-decomposition}
 For an edge \smash{$\xymatrix{t\ar@{-}[r]^k&t'}$} in $\mathbb{T}_n$ and $\varepsilon\in \{\pm 1\}$, we have
 \begin{align}\label{eq:FG-dec-gqca}
 \mu_{k;t}=\psi_{k;t}\bigl(\mathbf{\hat{Y}}_{k;t}^{\varepsilon}\bigr)^{\varepsilon}\circ\phi_{k;t;\varepsilon},
 \end{align}
 where $\phi_{k;t;\varepsilon}$ is the unique $\Z\big[q^{\pm\!\frac{1}{2}}\big]$-algebra isomorphism from $\mathcal{F}_{t'}$ to $\mathcal{F}_{\!t}$ taking $\mathbf{X}_{\!t'}(\alpha)$ to $\mathbf{X}_{t}\bigl(E_{\!k,\varepsilon}^{\tilde{B}_{t}R}\alpha\bigr)$ for any $\alpha\in\mathbb{Z}^{m}$.
\end{Proposition}
\begin{proof}
 By definition, we have
 \begin{align*}
 \mathbf{X}_{t'}(e_{k})
 &=\sum\limits^{r_{k}}_{s=0}h_{k;s}\bigl(q^{\frac{1}{2}}\bigr)\mathbf{X}_{t}(s[\varepsilon\mathbf{b}_{k;t}]_{+}+(r_{k}-s)[-\varepsilon\mathbf{b}_{k;t}]_{+}-e_{k})\\
 & =\sum\limits^{r_{k}}_{s=0}h_{k;s}\bigl(q^{\frac{1}{2}}\bigr)\mathbf{X}_{t}(s\varepsilon\mathbf{b}_{k;t}+r_{k}[-\varepsilon\mathbf{b}_{k;t}]_{+}-e_{k})\\
 &=\sum\limits^{r_{k}}_{s=0}h_{k;s}\bigl(q^{\frac{1}{2}}\bigr)\mathbf{X}_{t}\bigl(s\varepsilon\mathbf{b}_{k;t}+E_{k,\varepsilon}^{\tilde{B}_{t}R}e_{k}\bigr)\\
 &=\sum\limits^{r_{k}}_{s=0}h_{k;s}\bigl(q^{\frac{1}{2}}\bigr)q^{\frac{1}{2}\Lambda_{t}(s\varepsilon\mathbf{b}_{k;t},E_{k,\varepsilon}^{\tilde{B}_{t}R}e_{k})}\mathbf{X}_{t}\bigl(E^{\tilde{B}_{t}R}_{k,\varepsilon}e_{k}\bigr)\mathbf{X}_{t}(s\varepsilon\mathbf{b}_{k;t})\\
 &=\sum\limits^{r_{k}}_{s=0}h_{k;s}\bigl(q^{\frac{1}{2}}\bigr)q^{\frac{-s\varepsilon}{2d_{k}}}\mathbf{X}_{t}\bigl(E^{\tilde{B}_{t}R}_{k,\varepsilon}e_{k}\bigr)\mathbf{X}_{t}(s\varepsilon\mathbf{b}_{k;t})\\
 &=\mathbf{X}_{t}\bigl(E^{\tilde{B}_{t}R}_{k,\varepsilon}e_{k}\bigr)\sum\limits^{r_{k}}_{s=0}h_{k;s}\bigl(q^{\frac{1}{2}}\bigr)q^{\frac{-s\varepsilon}{2d_{k}}}\mathbf{X}_{t}(s\varepsilon\mathbf{b}_{k;t})\\
 &=\mathbf{X}_{t}\bigl(E^{\tilde{B}_{t}R}_{k,\varepsilon}e_{k}\bigr)\left(\sum\limits^{r_{k}}_{s=0}h_{k;s}\bigl(q^{\frac{1}{2}}\bigr)\bigl(q^{\frac{1}{2d_{k}}}\mathbf{X}_{t}(\varepsilon\mathbf{b}_{k;t})\bigr)^{s}\right)^{-\varepsilon\{-\varepsilon\}}\\
 &=\mathbf{X}_{t}\bigl(E^{\tilde{B}_{t}R}_{k,\varepsilon}e_{k}\bigr)\left(\sum\limits^{r_{k}}_{s=0}h_{k;s}\bigl(q^{\frac{1}{2}}\bigr)\bigl(q^{\frac{1}{2d_{k}}}\mathbf{\hat{Y}}_{k;t}^{\varepsilon}\bigr)^{s}\right)^{-\varepsilon\{-\varepsilon\}}.
 \end{align*}
 Since \smash{$\bigl(d_{k}\varepsilon f_k,\overline{E^{\tilde{B}_{t}R}_{k,\varepsilon}e_{i}}\bigr)_D=-\varepsilon\delta_{ik}$}, we have
 \begin{align*}
 \mathbf{X}_{t'}(e_{k})=\psi_{k;t}\bigl(\mathbf{\hat{Y}}_{k;t}^{\varepsilon}\bigr)^{\varepsilon}\bigl(\mathbf{X}_{t}\bigl(E^{\tilde{B}_{t}R}_{k,\varepsilon}e_{k}\bigr)\bigr).
 \end{align*}
 For $j\neq k$, we have
 \[
 \mathbf{X}_{t'}(e_{j})
 =\mathbf{X}_{t}(e_{j})
 =\mathbf{X}_{t}\bigl(E^{\tilde{B}_{t}R}_{k,\varepsilon}e_{j}\bigr)
 =\psi_{k;t}\bigl(\mathbf{\hat{Y}}_{k;t}^{\varepsilon}\bigr)^{\varepsilon}\bigl(\mathbf{X}_{t}\bigl(E^{\tilde{B}_{t}R}_{k,\varepsilon}e_{j}\bigr)\bigr).\tag*{\qed}
 \]\renewcommand{\qed}{}
\end{proof}

\begin{Remark}
When $R=I_n$, equation \eqref{eq:FG-dec-gqca} specializes to the Fock--Goncharov decomposition for quantum cluster algebras, where $\phi_{k;t;+}$ and $\phi_{k;t;-}$ are called the {\it tropical parts} of $\mu_{k;t}$ (cf.\ \cite[Section 3.3]{fock2009} and \cite[Section 6.3]{Keller2012}). By setting \smash{$q^{\frac12}=1$}, it further degenerate to the Fock--Goncharov decomposition for cluster algebras, see \cite[Section 2.3]{LMN2023}.
\end{Remark}

Fix a path \smash{$\xymatrix{t_{0}\ar@{-}[r]^{i_1}&t_{1}\ar@{-}[r]^{i_2}&t_{2}\ar@{-}[r]^{i_3}&\cdots\ar@{-}[r]^{i_k}&t_{k}}$} in $\mathbb{T}_n$ and denote it by $\mathbf{i}$. Let $\varepsilon_{j}$ be the common sign of components of $\mathbf{c}_{i_j;t_{j-1}}$ and $\mathbf{c}_{j}^{+}:=\varepsilon_{j}\mathbf{c}_{i_j;t_{j-1}}$ for $j\in[1,k]$.
\begin{Remark}
 The definition of $\varepsilon_j$ is opposite to the one in \cite{LMN2023}, where they define $\varepsilon_j$ to be the common sign of components of $\mathbf{c}_{i_j;t_j}$
 Our convention of $\varepsilon_{j}$ is the same as the one in \cite{Keller2012}, which is more convenient to study maximal green sequences and green-to-red sequences in cluster algebras \cite{KD2020}.
\end{Remark}
Now define
\begin{align*}
&\mu_{t_k}^{t_0}:=\mu_{i_1;t_0}\circ\mu_{i_2;t_1}\circ\cdots\mu_{i_k;t_{k-1}}\colon\ \mathcal{F}_{t_k}\to\mathcal{F}_{t_0},\\
&\psi(\mathbf{i}):=\psi_{i_1;t_0}\bigl(\mathbf{\hat{Y}}_{t_0}^{\mathbf{c}_{1}^+}\bigr)^{\varepsilon_1}\circ\psi_{i_2;t_0}\bigl(\mathbf{\hat{Y}}_{t_0}^{\mathbf{c}_{2}^+}\bigr)^{\varepsilon_2}\circ\cdots\circ\psi_{i_k;t_0}\bigl(\mathbf{\hat{Y}}_{t_0}^{\mathbf{c}_{k}^+}\bigr)^{\varepsilon_k}\colon\ \mathcal{F}_{t_0}\to\mathcal{F}_{t_0}.
\end{align*}
For each $j\in [1,k]$, we also set
\[
\phi_{t_{j}}^{t_0}:=\phi_{i_1;t_{0};\varepsilon_{1}}\circ\phi_{i_2;t_{1};\varepsilon_{2}}\circ\cdots\circ\phi_{i_j;t_{j-1};\varepsilon_{j}}\colon\ \mathcal{F}_{t_{j}}\to\mathcal{F}_{t_{0}}.
\]
It is straightforward to check that $\phi_{t_{j}}^{t_0}$ takes $\mathbf{X}_{t_{j}}(\alpha)$ to $\mathbf{X}_{t_{0}}\bigl(\tilde{G}_{t_{j}}\alpha\bigr)$ by \eqref{eq:g-matrix}, where $\alpha\in \Z^m$.
\begin{Lemma}\label{l:commutative-diag-isom}
 For each $j\in [1,k-1]$, we have
\[
 \phi_{t_{j}}^{t_{0}}\circ\psi_{i_{j+1};t_{j}}\bigl(\mathbf{\hat{Y}}^{\varepsilon_{j+1}}_{i_{j+1};t_{j}}\bigr)^{\varepsilon_{j+1}}
 =\psi_{i_{j+1};t_{0}}\bigl(\mathbf{\hat{Y}}_{t_{0}}^{\mathbf{c}_{j+1}^{+}}\bigr)^{\varepsilon_{j+1}}\circ\phi_{t_{j}}^{t_{0}}.
 \]
\end{Lemma}
\begin{proof}
It is equivalent to show
 \[
 \psi_{i_{j+1};t_{0}}\bigl(\mathbf{\hat{Y}}_{t_{0}}^{\mathbf{c}_{j+1}^{+}}\bigr)\circ\phi_{t_{j}}^{t_{0}}=\begin{cases}
 \phi_{t_{j}}^{t_{0}}\circ\psi_{i_{j+1};t_{j}}\bigl(\mathbf{\hat{Y}}_{i_{j+1};t_{j}}\bigr) & \text{if $\varepsilon_{j+1}=+1$},\\
 \phi_{t_{j}}^{t_{0}}\circ\psi_{i_{j+1};t_{j}}\bigl(\mathbf{\hat{Y}}^{-1}_{i_{j+1};t_{j}}\bigr) & \text{if $\varepsilon_{j+1}=-1$.}
 \end{cases}
 \]
 For any $\alpha=(\alpha_{1},\dots,\alpha_{m})^{\mathsf T}\in\mathbb{Z}^{m}$,
 \begin{gather*}
  \phi_{t_{j}}^{t_{0}}\circ\psi_{i_{j+1};t_{j}}\bigl(\mathbf{\hat{Y}}_{i_{j+1};t_{j}}^{\varepsilon_{j+1}}\bigr)(\mathbf{X}_{t_{j}}(\alpha))\\
  \qquad=\phi_{t_{j}}^{t_{0}}\left(\mathbf{X}_{t_{j}}(\alpha)\left(\sum\limits_{s=0}^{r_{i_{j+1}}} h_{i_{j+1},s}\bigl(q^{\frac{1}{2}}\bigr)\bigl(q^{\frac{1}{2d_{i_{j+1}}}}\mathbf{\hat{Y}}_{i_{j+1};t_{j}}^{\varepsilon_{j+1}}\bigr)^{s}\right)^{-\{d_{i_{j+1}}(\bar{\alpha},\varepsilon_{j+1}f_{i_{j+1}})_D\}}\right)\\
 \qquad =\mathbf{X}_{t_{0}}\bigl(\tilde{G}_{t_{j}}\alpha\bigr)\left(\sum\limits_{s=0}^{r_{i_{j+1}}} h_{i_{j+1},s}\bigl(q^{\frac{1}{2}}\bigr)\bigl(q^{\frac{1}{2d_{i_{j+1}}}}\mathbf{X}_{t_0}\bigl(\varepsilon_{j+1}\tilde{G}_{t_{j}}\tilde{B}_{t_{j}}f_{i_{j+1}}\bigr)\bigr)^{s}\right)^{-\{d_{i_{j+1}}(\bar{\alpha},\varepsilon_{j+1}f_{i_{j+1}})_D\}}\\
 \qquad =\mathbf{X}_{t_{0}}\bigl(\tilde{G}_{t_{j}}\alpha\bigr)\left(\sum\limits_{s=0}^{r_{i_{j+1}}} h_{i_{j+1},s}\bigl(q^{\frac{1}{2}}\bigr)\bigl(q^{\frac{1}{2d_{i_{j+1}}}}\mathbf{X}_{t_0}\bigl(\varepsilon_{j+1}\tilde{B}_{t_{0}}C_{t_{j}}f_{i_{j+1}}\bigr)\bigr)^{s}\right)^{-\{d_{i_{j+1}}(\bar{\alpha},\varepsilon_{j+1}f_{i_{j+1}})_D\}}\\
 \qquad =\mathbf{X}_{t_{0}}\bigl(\tilde{G}_{t_{j}}\alpha\bigr)\left(\sum\limits_{s=0}^{r_{i_{j+1}}} h_{i_{j+1},s}\bigl(q^{\frac{1}{2}}\bigr)\bigl(q^{\frac{1}{2d_{i_{j+1}}}}\mathbf{X}_{t_0}\bigl(\tilde{B}_{t_{0}}\mathbf{c}_{j+1}^{+}\bigr)\bigr)^{s}\right)^{-\{d_{i_{j+1}}(\bar{\alpha},\varepsilon_{j+1}f_{i_{j+1}})_D\}}\\
 \qquad =\mathbf{X}_{t_{0}}\bigl(\tilde{G}_{t_{j}}\alpha\bigr)\left(\sum\limits_{s=0}^{r_{i_{j+1}}} h_{i_{j+1},s}\bigl(q^{\frac{1}{2}}\bigr)\bigl(q^{\frac{1}{2d_{i_{j+1}}}}\mathbf{\hat{Y}}_{t_{0}}^{\mathbf{c}_{j+1}^{+}}\bigr)^{s}\right)^{-\{\alpha_{i_{j+1}}\varepsilon_{j+1}\}},
 \end{gather*}
 where the third equality follows from \eqref{duality_quantum_G_C_matrix1}.
 On the other hand,
 \begin{gather*}
 \psi_{i_{j+1};t_{0}}\bigl(\mathbf{\hat{Y}}_{t_{0}}^{\mathbf{c}_{j+1}^{+}}\bigr)\circ\phi_{t_{j}}^{t_{0}}(\mathbf{X}_{t_{j}}(\alpha))\\
\qquad =\psi_{i_{j+1};t_{0}}\bigl(\mathbf{\hat{Y}}_{t_{0}}^{\mathbf{c}_{j+1}^{+}}\bigr)\bigl(\mathbf{X}_{t_{0}}\bigl(\tilde{G}_{t_{j}}\alpha\bigr)\bigr)\\
\qquad =\mathbf{X}_{t_{0}}\bigl(\tilde{G}_{t_{j}}\alpha\bigr)\left(\sum\limits_{s=0}^{r_{i_{j+1}}} h_{i_{j+1},s}\bigl(q^{\frac{1}{2}}\bigr)\bigl(q^{\frac{1}{2d_{i_{j+1}}}}\mathbf{\hat{Y}}_{t_{0}}^{\mathbf{c}_{j+1}^{+}}\bigr)^{s}\right)^{-\{d_{i_{j+1}}(\overline{\tilde{G}_{t_{j}}\alpha},\mathbf{c}_{j+1}^{+})_D\}}\\
\qquad =\mathbf{X}_{t_{0}}\bigl(\tilde{G}_{t_{j}}\alpha\bigr)\left(\sum\limits_{s=0}^{r_{i_{j+1}}} h_{i_{j+1},s}\bigl(q^{\frac{1}{2}}\bigr)\bigl(q^{\frac{1}{2d_{i_{j+1}}}}\mathbf{\hat{Y}}_{t_{0}}^{\mathbf{c}_{j+1}^{+}}\bigr)^{s}\right)^{-\{\alpha_{i_{j+1}}\varepsilon_{j+1}\}}\qquad\text{(by \eqref{D_form_duality})}.
 \end{gather*}
 This completes the proof.
\end{proof}

\begin{Proposition}\label{long-quantum-decomposition}
 Keep the notation as above, we have
 $\mu_{t_k}^{t_0}=\psi(\mathbf{i})\circ\phi_{t_{k}}^{t_0}\colon  \mathcal{F}_{t_{k}}\to\mathcal{F}_{t_{0}}$.
\end{Proposition}
\begin{proof}
 By using Proposition \ref{quantum-decomposition} and Lemma \ref{l:commutative-diag-isom}, we have
 \begin{align*}
 \mu_{t_k}^{t_0}
 ={}&\mu_{i_1;t_0}\circ\mu_{i_2;t_1}\circ\cdots\circ\mu_{i_k;t_{k-1}}\\
 ={}&\psi_{i_{1};t_{0}}\bigl(\mathbf{\hat{Y}}_{i_{1};t_{0}}^{\varepsilon_{1}}\bigr)^{\varepsilon_{1}}\circ\phi_{i_{1};t_{0};\varepsilon_{1}}\circ\psi_{i_{2};t_{1}}
 \bigl(\mathbf{\hat{Y}}_{i_{2};t_{1}}^{\varepsilon_{2}}\bigr)^{\varepsilon_{2}}\circ\phi_{i_{2};t_{1};\varepsilon_{2}}\\
 &\circ\cdots\circ\psi_{i_{k};t_{k-1}}
 \bigl(\mathbf{\hat{Y}}_{i_{k};t_{k-1}}^{\varepsilon_{k}}\bigr)^{\varepsilon_{k}}\circ\phi_{i_{k};t_{k-1};\varepsilon_{k}}\\
 ={}&\psi_{i_{1};t_{0}}\bigl(\mathbf{\hat{Y}}_{t_{0}}^{\mathbf{c}_{1}^{+}}\bigr)^{\varepsilon_{1}}\circ\psi_{i_{2};t_{1}}\bigl(\mathbf{\hat{Y}}_{t_{0}}^{\mathbf{c}_{2}^{+}}\bigr)^{\varepsilon_{2}}\circ\phi_{i_{1};t_{0};\varepsilon_{1}}\circ\phi_{i_{2};t_{1};\varepsilon_{2}}\\
 &\circ\cdots\circ\psi_{i_{k};t_{k-1}}(\mathbf{\hat{Y}}_{i_{k};t_{k-1}}^{\varepsilon_{k}})^{\varepsilon_{k}}\circ\phi_{i_{k};t_{k-1};\varepsilon_{k}}\\
 &\vdots\\
 ={}&\psi_{i_{1};t_{0}}\bigl(\mathbf{\hat{Y}}_{t_{0}}^{\mathbf{c}_{1}^{+}}\bigr)^{\varepsilon_{1}}\circ\psi_{i_{2};t_{1}}\bigl(\mathbf{\hat{Y}}_{t_{0}}^{\mathbf{c}_{2}^{+}}\bigr)^{\varepsilon_{2}}\circ\cdots\circ\psi_{i_{k};t_{0}}
 \bigl(\mathbf{\hat{Y}}_{t_{0}}^{\mathbf{c}_{k}^{+}}\bigr)^{\varepsilon_{k}}\circ\phi_{t_{k}}^{t_0}\\
 ={}&\psi(\mathbf{i})\circ\phi_{t_{k}}^{t_0}.\tag*{\qed}
 \end{align*}\renewcommand{\qed}{}
\end{proof}

\subsection[F-polynomials]{$\boldsymbol{F}$-polynomials}\label{s:F-poly_gen_quantum_cluster_alg}

 In this section, we define $F$-polynomials for the generalized quantum cluster algebra $\mathcal{A}_q(\Sigma_{t_0})$ and prove their polynomial property under a mild condition.
Recall that for $i,j\in[1,n]$, we have
\begin{align}\label{eq:quasi-com}
\mathbf{\hat{Y}}_{i,t_{0}}\mathbf{\hat{Y}}_{j,t_{0}}=q^{d_{i}^{-1}b_{ij}}\mathbf{\hat{Y}}_{j,t_{0}}\mathbf{\hat{Y}}_{i,t_{0}}.
\end{align}
The quasi-commutative relations \eqref{eq:quasi-com} depend only on the entries of $D$ and the principal part~$B$ of $\tilde{B}_{t_0}$.
Let $\mathcal{T}_{DB}$ be the quantum torus associated to $DB$ whose underlying space is the free~\smash{$\mathbb{Z}\big[q^{\frac{1}{2}}\big]$}-module with basis $\{\mathbf{Z}(\alpha)| \alpha\in\mathbb{Z}^{n}\}$, and its multiplication is given by
\[\mathbf{Z}(\alpha)\mathbf{Z}(\beta)=q^{\frac{1}{2}\alpha^{\mathsf T}DB\beta}\mathbf{Z}(\alpha+\beta).\]
Denote by $\mathcal{F}_{DB}$ the fraction skew field of $\mathcal{T}_{DB}$ and $\mathbf{Z}_{i}=\mathbf{Z}(f_{i})$.

Recall that \smash{$\xymatrix{t_{0}\ar@{-}[r]^{i_1}&t_{1}\ar@{-}[r]^{i_2}&t_{2}\ar@{-}[r]^{i_3}&\cdots\ar@{-}[r]^{\!\!i_k}&t_{k}}$} is a path in $\mathbb{T}_n$, and $\varepsilon_{j}$ is the common sign of components of $\mathbf{c}_{i_j;t_{j-1}}$. For simplicity of notation, we also denote by
\begin{alignat*}{6}
& d_{(j)}=d_{i_j},\qquad&&
r_{(j)}=r_{i_j},  \qquad &&  \mathbf{c}_{j}=\mathbf{c}_{i_j;t_{j-1}},\qquad&&    \mathbf{c}_{j}^+=\varepsilon_{j}\mathbf{c}_{j}, \qquad&&   \mathbf{\hat{c}}_{j}^+=B\mathbf{c}_{j}^+,& \\
& \mathbf{\tilde{g}}_{j}=\mathbf{\tilde{g}}_{i_j;t_j}, \qquad&&   \mathbf{g}_{j}=\overline{\mathbf{\tilde{g}}_{j}}.&&&&&&&
\end{alignat*}

We first define a set of elements $\{ L_{i,j}\mid i,j\in [1,k]\}$ of $\mathcal{F}_q$ by the initial condition
\[L_{1,i}:=\mathbf{\hat{Y}}_{t_{0}}^{\mathbf{c}_{i}^{+}}\left(\sum\limits^{r_{(1)}}_{s=0}h_{i_{1},s}
\bigl(q^{\frac{1}{2}}\bigr)\bigl(q^{\frac{1}{2d_{(1)}}}\mathbf{\hat{Y}}_{t_0}^{\mathbf{c}_{1}^{+}}\bigr)^{s}\right)^{-\varepsilon_{1}\{(d_{(1)}\mathbf{c}_{1}^{+},\mathbf{\hat{c}}_{i}^{+})_{D}\}} \qquad \text{for}\  i\in [1,k]
\]
with recurrence relations
\[
 L_{j+1,i}=L_{j,i}\left(\sum\limits^{r_{(j+1)}}_{s=0}h_{i_{j+1},s}\bigl(q^{\frac{1}{2}}\bigr)
 \bigl(q^{\frac{1}{2d_{(j+1)}}}L_{j,j+1}\bigr)^{s} \right)^{-\varepsilon_{j+1}\{(d_{(j+1)}\mathbf{c}_{j+1}^{+},\mathbf{\hat{c}}_{i}^{+})_{D}\}} \quad \text{for}\  j \in [1,k-1].
\]
Then set
\begin{align*}
 &L_1=\sum\limits^{r_{(1)}}_{s=0}h_{i_{1},s}\bigl(q^{\frac{1}{2}}\bigr)\bigl(q^{\frac{1}{2d_{(1)}}}\mathbf{\hat{Y}}_{t_{0}}^{\mathbf{c}_{1}^{+}}\bigr)^{s},\\
 &L_{l+1}=\sum\limits_{s=0}^{r_{(l+1)}}h_{i_{l+1},s}\bigl(q^{\frac{1}{2}}\bigr)\bigl(q^{\frac{1}{2d_{(l+1)}}}L_{l,l+1}\bigr)^{s},\qquad l\in[1,k-1].
 \end{align*}

\begin{Lemma}\label{product_of_F_polynomial}
 Keep the notation as above. We have
 \begin{align}\label{eq:f-poly-well-def}
 \mathbf{X}_{t_{0}}(-\mathbf{\tilde{g}}_{k})\mu_{t_{k}}^{t_{0}}(\mathbf{X}_{i_{k};t_{k}})=\prod\limits^{\longrightarrow}_{j\in[1,k]}L_{j}^{-\varepsilon_{j}\{d_{(j)}(\mathbf{c}_{j}^{+},\mathbf{g}_{k})_D\}}.
 \end{align}

\end{Lemma}
\begin{proof}
 Let $\mathbf{i}_{j}$ be the sub-sequence \smash{$\xymatrix{t_{0}\ar@{-}[r]^{i_1}&t_{1}\ar@{-}[r]^{i_2}&t_{2}\ar@{-}[r]^{i_3}&\cdots\ar@{-}[r]^{i_j}&t_{j}}$} for $j\in [1,k]$.
 By Proposition \ref{long-quantum-decomposition},
 \begin{align*}
 \mathbf{X}_{t_{0}}(-\mathbf{\tilde{g}}_{k})\mu_{t_{k}}^{t_{0}}(\mathbf{X}_{i_{k};t_{k}})
 ={}&\mathbf{X}_{t_{0}}(-\mathbf{\tilde{g}}_{k})\psi(\mathbf{i})\circ\phi_{t_{k}}^{t_0}(\mathbf{X}_{i_{k};t_{k}})\\
 ={}&\mathbf{X}_{t_{0}}(-\mathbf{\tilde{g}}_{k})\psi(\mathbf{i})(\mathbf{X}_{t_{0}}(\mathbf{\tilde{g}}_{k}))\\
 ={}&\mathbf{X}_{t_{0}}(-\mathbf{\tilde{g}}_{k})\psi(\mathbf{i}_{1})(\mathbf{X}_{t_{0}}(\mathbf{\tilde{g}}_{k}))\psi(\mathbf{i}_{1})(\mathbf{X}_{t_{0}}(-\mathbf{\tilde{g}}_{k}))\psi(\mathbf{i}_{2})(\mathbf{X}_{t_{0}}(\mathbf{\tilde{g}}_{k}))\\
 &\cdots\psi(\mathbf{i}_{k-1})(\mathbf{X}_{t_{0}}(-\mathbf{\tilde{g}}_{k}))\psi(\mathbf{i}_{k})(\mathbf{X}_{t_{0}}(\mathbf{\tilde{g}}_{k})).
 \end{align*}
 We first claim that $L_{j,i}=\psi(\mathbf{i}_{j})\bigl(\mathbf{\hat{Y}}_{t_{0}}^{\mathbf{c}_{i}^{+}}\bigr)$
 for $i,j\in[1,k]$. We prove it by induction on $j$.
 For $j=1$, we have
 \begin{align*}
 \psi(\mathbf{i}_{1})\bigl(\mathbf{\hat{Y}}_{t_{0}}^{\mathbf{c}_{i}^{+}}\bigr)
 &=\psi_{i_{1}:t_{0}}\bigl(\mathbf{\hat{Y}}_{t_{0}}^{\mathbf{c}_{1}^{+}}\bigr)^{\varepsilon_{1}}\bigl(\mathbf{\hat{Y}}_{t_{0}}^{\mathbf{c}_{i}^{+}}\bigr)=\mathbf{\hat{Y}}_{t_{0}}^{\mathbf{c}_{i}^{+}}\left(\sum\limits^{r_{(1)}}_{s=0}h_{i_{1},s}\bigl(q^{\frac{1}{2}}\bigr)
 \bigl(q^{\frac{1}{2d_{(1)}}}\mathbf{\hat{Y}}_{t_{0}}^{\mathbf{c}_{1}^{+}}\bigr)^{s}\right)^{-\varepsilon_{1}\{(\mathbf{\hat{c}}_{i}^{+},d_{(1)}\mathbf{c}_{1}^{+})_{D}\}}\\
 &=L_{1,i}.
 \end{align*}
 Suppose that $L_{j,i}=\psi(\mathbf{i}_{j})\bigl(\mathbf{\hat{Y}}_{t_{0}}^{\mathbf{c}_{i}^{+}}\bigr)$ for any $i\in[1,k]$, then
 \begin{align*}
 \psi(\mathbf{i}_{j+1})\bigl(\mathbf{\hat{Y}}_{t_{0}}^{\mathbf{c}_{i}^{+}}\bigr)
 ={}&\psi(\mathbf{i}_{j})\bigl(\psi_{i_{j+1};t_{0}}\bigl(\mathbf{\hat{Y}}_{t_{0}}^{\mathbf{c}_{j+1}^{+}}\bigr)^{\varepsilon_{j+1}}\bigl(\mathbf{\hat{Y}}_{t_{0}}^{\mathbf{c}_{i}^{+}}\bigr)\bigr)\\
 ={}&\psi(\mathbf{i}_{j})\left(\bigl(\mathbf{\hat{Y}}_{t_{0}}^{\mathbf{c}_{i}^{+}}\bigr)\left(\sum\limits^{r_{(j+1)}}_{s=0}h_{i_{j+1},s}
 \bigl(q^{\frac{1}{2}}\bigr)\bigl(q^{\frac{1}{2d_{(j+1)}}}\mathbf{\hat{Y}}_{t_{0}}^{\mathbf{c}_{j+1}^{+}}\bigr)^{s}\right)^{-\varepsilon_{j+1}
 \{(\mathbf{\hat{c}}_{i}^{+},d_{(j+1)}\mathbf{c}_{j+1}^{+})_{D}\}}\right)\\
 ={}&\psi(\mathbf{i}_{j})\bigl(\mathbf{\hat{Y}}_{t_{0}}^{\mathbf{c}_{i}^{+}}\bigr)\\
 &\times\left(\sum\limits^{r_{(j+1)}}_{s=0}h_{i_{j+1},s}\bigl(q^{\frac{1}{2}}\bigr)
 \bigl(q^{\frac{1}{2d_{(j+1)}}}\psi(\mathbf{i}_{j})\bigl(\mathbf{\hat{Y}}_{t_{0}}^{\mathbf{c}_{j+1}^{+}}\bigr)\bigr)^{s}\right)^{-\varepsilon_{j+1}\{(\mathbf{\hat{c}}_{i}^{+},d_{(j+1)}\mathbf{c}_{j+1}^{+})_{D}\}}\\
 ={}&L_{j,i}\left(\sum\limits^{r_{(j+1)}}_{s=0}h_{i_{j+1},s}\bigl(q^{\frac{1}{2}}\bigr)\bigl(q^{\frac{1}{2d_{(j+1)}}}L_{j,j+1}\bigr)^{s}\right)^{-\varepsilon_{j+1}
 \{(\mathbf{\hat{c}}_{i}^{+},d_{(j+1)}\mathbf{c}_{j+1}^{+})_{D}\}}\\
 ={}&L_{j+1,i}.
 \end{align*}
 This completes the proof of the claim.
 A direct computation shows that
 \begin{align*}
 \mathbf{X}_{t_0}(-\mathbf{\tilde{g}}_{k})\psi(\mathbf{i}_{1})(\mathbf{X}_{t_0}(\mathbf{\tilde{g}}_{k}))
 &=\mathbf{X}_{t_0}(-\mathbf{\tilde{g}}_{k})\psi_{i_{1};t_{0}}\bigl(\mathbf{\hat{Y}}_{t_{0}}^{\mathbf{c}_{1}^{+}}\bigr)^{\varepsilon_{1}}(\mathbf{X}_{t_0}(\mathbf{\tilde{g}}_{k}))\\
 &=\left(\sum\limits^{r_{(1)}}_{s=0}h_{i_{1},s}\bigl(q^{\frac{1}{2}}\bigr)\bigl(q^{\frac{1}{2d_{(1)}}}\mathbf{\hat{Y}}_{t_{0}}^{\mathbf{c}_{1}^{+}}\bigr)^{s}\right)^{-\varepsilon_{1}\{(\mathbf{g}_{k},d_{(1)}\mathbf{c}_{1}^{+})_{D}\}}\\
 &=L_{1}^{-\varepsilon_{1}\{(\mathbf{g}_{k},d_{(1)}\mathbf{c}_{1}^{+})_{D}\}}.
 \end{align*}
 For $j\in[1,k-1]$, we have
 \begin{gather*}
  \psi(\mathbf{i}_{j})(\mathbf{X}_{t_{0}}(-\mathbf{\tilde{g}}_{k}))\psi(\mathbf{i}_{j+1})(\mathbf{X}_{t_{0}}(\mathbf{\tilde{g}}_{k}))\\
  \qquad=\psi(\mathbf{i}_{j})\bigl(\mathbf{X}_{t_{0}}(-\mathbf{\tilde{g}}_{k})\psi_{i_{j+1};t_{j}}\bigl(\mathbf{\hat{Y}}_{t_0}^{\mathbf{c}_{j+1}^+}\bigr)^{\varepsilon_{j+1}}
  (\mathbf{X}_{t_{0}}(\mathbf{\tilde{g}}_{k}))\bigr)\\
 \qquad  =\psi(\mathbf{i}_{j})\left(\left(\sum\limits^{r_{(j+1)}}_{s=0}h_{i_{j+1},s}\bigl(q^{\frac{1}{2}}\bigr)
 \bigl(q^{\frac{1}{2d_{(j+1)}}}\mathbf{\hat{Y}}_{t_{0}}^{\mathbf{c}_{j+1}^{+}}\bigr)^{s}\right)
 ^{-\varepsilon_{j+1}\{(\mathbf{g}_{k},d_{(j+1)}\mathbf{c}_{j+1}^{+})_{D}\}}\right)\\
\qquad  =\left(\sum\limits^{r_{(j+1)}}_{s=0}h_{i_{j+1},s}\bigl(q^{\frac{1}{2}}\bigr)\bigl(q^{\frac{1}{2d_{(j+1)}}}\psi(\mathbf{i}_{j})\bigl(\mathbf{\hat{Y}}_{t_{0}}^{\mathbf{c}_{j+1}^{+}}\bigr)\bigr)^{s}\right)^{-\varepsilon_{j+1}\{(\mathbf{g}_{k},d_{(j+1)}\mathbf{c}_{j+1}^{+})_{D}\}}\\
\qquad  =\left(\sum\limits^{r_{(j+1)}}_{s=0}h_{i_{j+1},s}\bigl(q^{\frac{1}{2}}\bigr)\bigl(q^{\frac{1}{2d_{(j+1)}}}L_{j,j+1}\bigr)^{s}\right)^{-\varepsilon_{j+1}\{(\mathbf{g}_{k},d_{(j+1)}\mathbf{c}_{j+1}^{+})_{D}\}}\\
\qquad  =L_{j+1}^{-\varepsilon_{j+1}\{(\mathbf{g}_{k},d_{(j+1)}\mathbf{c}_{j+1}^{+})_{D}\}}.
 \end{gather*}
 It follows that
 \begin{align*}
 \mathbf{X}_{t_{0}}(-\mathbf{\tilde{g}}_{k})\mu_{t_{k}}^{t_{0}}(\mathbf{X}_{i_{k};t_{k}})
 ={}&\mathbf{X}_{t_{0}}(-\mathbf{\tilde{g}}_{k})\psi(\mathbf{i}_{1})(\mathbf{X}_{t_{0}}(\mathbf{\tilde{g}}_{k}))\psi(\mathbf{i}_{1})(\mathbf{X}_{t_{0}}(-\mathbf{\tilde{g}}_{k}))\psi(\mathbf{i}_{2})(\mathbf{X}_{t_{0}}(\mathbf{\tilde{g}}_{k}))\\
 &\cdots\psi(\mathbf{i}_{k-1})(\mathbf{X}_{t_{0}}(-\mathbf{\tilde{g}}_{k}))\psi(\mathbf{i}_{k})(\mathbf{X}_{t_{0}}(\mathbf{\tilde{g}}_{k}))\\
 ={}&L_{1}^{-\varepsilon_{1}\{(\mathbf{g}_{k},d_{(1)}\mathbf{c}_{1}^{+})_{D}\}}L_{2}^{-\varepsilon_{2}\{(\mathbf{g}_{k},d_{(2)}\mathbf{c}_{2}^{+})_{D}\}}\cdots L_{k}^{-\varepsilon_{k}\{(\mathbf{g}_{k},d_{(k)}\mathbf{c}_{k}^{+})_{D}\}}.\tag*{\qed}
 \end{align*}\renewcommand{\qed}{}
\end{proof}

\begin{Definition}
 There is a unique element $F_{i_k;t_k}:=F_{i_k;t_k}[\mathbf{Z}_1,\dots, \mathbf{Z}_n]$ in the skew field $\mathcal{F}_{DB}$ of~$\mathcal{T}_{DB}$ such that
 \begin{align}\label{eq:gupta-formula}
 \prod\limits^{\longrightarrow}_{j\in[1,k]}L_{j}^{-\varepsilon_{j}\{d_{(j)}(\mathbf{c}_{j}^{+},\mathbf{g}_{k})_D\}}=F_{i_k;t_k}\big[\hat{\mathbf{Y}}_{1;t_0},\dots, \hat{\mathbf{Y}}_{n;t_0}\big].
 \end{align}
 The element $F_{i_k,t_k}$ is called the {\it $F$-polynomial} of $\mathbf{X}_{t_k}(e_{i_k})$ whenever $F_{i_k;t_k}$ is a polynomial in $\mathbf{Z}_1,\dots,\mathbf{Z}_n$. With the help of $F_{i_k;t_k}$, we can rewrite \eqref{eq:f-poly-well-def} as
 \[
 \mathbf{X}_{i_k;t_k}=\mathbf{X}_{t_0}(\tilde{\mathbf{g}}_k)F_{i_k;t_k}\big[\hat{\mathbf{Y}}_{1;t_0},\dots, \hat{\mathbf{Y}}_{n;t_0}\big].
 \]

\end{Definition}
\begin{Remark}
 Recall that the $(\mathbf{r},\mathbf{z})$-$C$-pattern and $(\mathbf{r},\mathbf{z})$-$G$-pattern are uniquely determined by $B$, $R$ and $t_0\in \mathbb{T}_n$. It follows that the element $F_{i_k;t_k}$ only depends on $B$, $R$, $D$ and $t_0$.
\end{Remark}

The following is the main result of this section.
\begin{Theorem}\label{t:F-polynomial}\qquad
\begin{itemize}\itemsep=0pt
 \item[$(1)$] The element $F_{i_k;t_k}$ is a Laurent polynomial in $\mathbf{Z}_1,\dots, \mathbf{Z}_n$.
 \item[$(2)$] Suppose that $h_{i,s}(1)>0$ for each $i\in [1,n]$ and $s\in [1,r_i-1]$, then $F_{i_k;t_k}$ is a polynomial in $\mathbf{Z}_1,\dots, \mathbf{Z}_n$.
\end{itemize}

\end{Theorem}
\begin{proof}
 Since $F_{i_k;t_k}$ only depends on $B$, $D$ and $R$, but not on coefficients, it suffices to prove the statement for a particular choice of coefficients. Let \smash{$\tilde{B}^\bullet=\big[\begin{smallmatrix}
 B\\ I_n
 \end{smallmatrix}\big]$} be the $2n\times n$ matrix, where~$I_n$ is the identity matrix. There is a $(2n\times 2n)$-skew symmetric matrix $\Lambda^{\bullet}$ such that $\bigl(\tilde{B}^{\bullet},\Lambda^{\bullet}\bigr)$ is a compatible pair with \smash{$\bigl(\tilde{B}^{\bullet}\bigr)^{\mathsf T}\Lambda^{\bullet}=[D\quad 0]$} (cf.\ \cite[Example 0.5]{Berenstein_Zelevinsky2005}). It is clear that the \smash{$\mathbb{Z}\big[q^{\pm\frac{1}{2}}\big]$}-algebra generated by $\bigl\{\mathbf{\hat{Y}}_{t_0}^{\alpha}:=\mathbf{X}_{t_{0}}\bigl(\tilde{B}^{\bullet}\alpha\bigr)\mid \alpha\in\mathbb{Z}^{n}\bigr\}$ is isomorphic to $\mathcal{T}_{DB}$.

 According to Lemma \ref{product_of_F_polynomial}, there are two polynomials $A(\mathbf{Z}_{1},\dots,\mathbf{Z}_{n}),P(\mathbf{Z}_{1},\dots,\mathbf{Z}_{n})\in \mathcal{T}_{DB}$ with coefficients in \smash{$\mathbb{N}\big[q^{\pm\frac{1}{2}}\big]$} for \smash{$h_{i,s}\bigl(q^{\frac{1}{2}}\bigr)$}\footnote{Here we temporary regard $h_{i,s}\bigl(q^{\frac{1}{2}}\bigr)$ as a variable.} and $\mathbf{Z}_{1},\dots,\mathbf{Z}_{n}$ such that
 \begin{align}\label{fraction_of_F-polynomial}
 F_{i_{k};t_{k}}\bigl(\mathbf{\hat{Y}}_{1;t_{0}},\dots,\mathbf{\hat{Y}}_{n;t_{0}}\bigr)=A\bigl(\mathbf{\hat{Y}}_{1;t_{0}},\dots,\mathbf{\hat{Y}}_{n;t_{0}}\bigr)P\bigl(\mathbf{\hat{Y}}_{1;t_{0}},\dots,\mathbf{\hat{Y}}_{n;t_{0}}\bigr)^{-1}.
 \end{align}

 By \cite[Theorem 3.1]{BCDX2023}, $F_{i_{k};t_{k}}\bigl(\mathbf{\hat{Y}}_{1;t_{0}},\dots,\mathbf{\hat{Y}}_{n;t_{0}}\bigr)$ is also a Laurent polynomial in $\mathbf{X}_{1;t_{0}},\dots,\mathbf{X}_{2n;t_{0}}$. Hence all of \smash{$F_{i_{k};t_{k}}\bigl(\mathbf{\hat{Y}}_{1;t_{0}},\dots,\mathbf{\hat{Y}}_{n;t_{0}}\bigr)$}, \smash{$A\bigl(\mathbf{\hat{Y}}_{1;t_{0}},\dots,\mathbf{\hat{Y}}_{n;t_{0}}\bigr)$} and \smash{$P\bigl(\mathbf{\hat{Y}}_{1;t_{0}},\dots,\mathbf{\hat{Y}}_{n;t_{0}}\bigr)$} are Laurent polynomials in $\mathbf{X}_{1;t_{0}},\dots,\mathbf{X}_{2n;t_{0}}$. Taking the Newton polytopes of both sides of \eqref{fraction_of_F-polynomial} as the Laurent polynomials in $\mathbf{X}_{1;t_0},\dots, \mathbf{X}_{2n;t_0}$, we obtain
$\text{New}(F_{i_{k};t_{k}})+\text{New}(P)=\text{New}(A)$,
 where the sum is the Minkowski sum. By the definition of Minkowski sum, we conclude that every exponent vector of $F_{i_{k};t_{k}}$ is a $\mathbb{Q}$-linear combination of exponents vectors of $A$ and $P$.
 In particular, the exponent vectors of
 \smash{$F_{i_{k};t_{k}}\bigl(\mathbf{\hat{Y}}_{1;t_{0}},\dots,\mathbf{\hat{Y}}_{n;t_{0}}\bigr)$} as a Laurent polynomial in $\mathbf{X}_{1;t_0}$,$\dots$, $\mathbf{X}_{2n;t_0}$ can be expressed as $\mathbb{Q}$-linear combinations of $\tilde{B}^{\bullet}f_{1},\dots,\tilde{B}^{\bullet}f_{n}$. However, the last $n$ coordinates of the vectors $\tilde{B}^{\bullet}f_{i}$ is given by the standard basis $f_1,\dots, f_{n}\in\mathbb{Z}^{n}$, so each exponent vector must be a $\mathbb{Z}$-linear combination of $\tilde{B}^{\bullet}f_{1},\dots,\tilde{B}^{\bullet}f_{n}$. This prove that $F_{i_{k};t_{k}}\bigl(\mathbf{\hat{Y}}_{1;t_{0}},\dots,\mathbf{\hat{Y}}_{n;t_{0}}\bigr)$ is a Laurent polynomial in $\mathbf{\hat{Y}}_{1;t_{0}},\dots,\mathbf{\hat{Y}}_{n;t_{0}}$. Hence, $F_{i_k;t_k}[\mathbf{Z}_1,\dots, \mathbf{Z}_n]$ is a Laurent polynomial in $\mathbf{Z}_{1},\dots,\mathbf{Z}_{n}$.

 Since $A(\mathbf{Z}_{1},\dots,\mathbf{Z}_{n}),P(\mathbf{Z}_{1},\dots,\mathbf{Z}_{n})$ have coefficients in $\mathbb{N}\big[q^{\pm\frac{1}{2}}\big]$ for $h_{i,s}\bigl(q^{\frac{1}{2}}\bigr)$ and $\mathbf{Z}_{1},\dots,\mathbf{Z}_{n}$. By the assumption $h_{i,s}(1)>0$ for $i\in[1,n]$ and $s\in[1,r_{i}-1]$, by setting \smash{$q^{\frac{1}{2}}=1$} does not shrink $\text{New}\bigl(F_{i_{k};t_{k}}\bigl(\mathbf{\hat{Y}}_{1;t_{0}},\dots,\mathbf{\hat{Y}}_{n;t_{0}}\bigr)\bigr)$. That is
 \[
 \text{New}\bigl(F_{i_{k};t_{k}}\bigl(\mathbf{\hat{Y}}_{1;t_{0}},\dots,\mathbf{\hat{Y}}_{n;t_{0}}\bigr)\bigr)
 =\text{New}\bigl(F_{i_{k};t_{k}}\bigl(\mathbf{\hat{Y}}_{1;t_{0}},\dots,\mathbf{\hat{Y}}_{n;t_{0}}\bigr)|_{q^{\frac{1}{2}}=1}\bigr).\]
By Proposition \ref{p:separation-formula-gca} and Lemma \ref{l:g-vector-comparison}, we conclude that \[F_{i_{k};t_{k}}\bigl(\mathbf{\hat{Y}}_{1;t_{0}},\dots,\mathbf{\hat{Y}}_{n;t_{0}}\bigr)|_{q^{\frac{1}{2}}=1}=F_{i_k;t_k}(\mathbf{\hat{y}},\mathbf{z})|_{z_{i,s}=h_{i,s}(1),i\in [1,n],s\in [1,r_i-1]},\] where $F_{i_k;t_k}(\mathbf{y},\mathbf{z})$ is the $F$-polynomial of the cluster variable $x_{i_k;t_k}$ of the corresponding generalized cluster algebra $\mathcal{A}^\bullet$ with principal coefficients.
 Thus $\text{New}\bigl(F_{i_{k};t_{k}}\bigl(\mathbf{\hat{Y}}_{1;t_{0}},\dots,\mathbf{\hat{Y}}_{n;t_{0}}\bigr)\bigr)$ does not contain any points with negative coordinates. It follows that \smash{$F_{i_{k};t_{k}}\bigl(\mathbf{\hat{Y}}_{1;t_{0}},\dots,\mathbf{\hat{Y}}_{n;t_{0}}\bigr)$} is a~polynomial in $\mathbf{\hat{Y}}_{1;t_{0}},\dots,\mathbf{\hat{Y}}_{n;t_{0}}$. Hence, $F_{i_{k};t_{k}}[\mathbf{Z}_1,\dots,\mathbf{Z}_{n}]$ is a polynomial in $\mathbf{Z}_{1},\dots,\mathbf{Z}_{n}$.
\end{proof}

As direct consequences of Theorem \ref{t:F-polynomial} and Lemma \ref{product_of_F_polynomial}, we have the following.

\begin{Corollary}[separation formula]\label{c:separation-formula}
Suppose that $h_{j,s}(1)>0$ for each $j\in [1,n]$ and $s\in [1,r_j-1]$. For each $i\in [1,n]$ and $t\in \mathbb{T}_n$, let $F_{i;t}[Z_1,\dots, Z_n]$ be the associated $F$-polynomial of~$\mathbf{X}_{i;t}$ and $\tilde{\mathbf{g}}_{i;t}$ the $g$-vector of $\mathbf{X}_{i;t}$. We have
\[\mathbf{X}_{i;t}=\mathbf{X}_{t_0}(\mathbf{\tilde{g}}_{i;t})F_{i;t}\bigl(\mathbf{\hat{Y}}_{1;t_{0}},\dots,\mathbf{\hat{Y}}_{n;t_{0}}\bigr).\]
\end{Corollary}
\begin{Corollary}\label{c:constant-term-1}
 Suppose that $h_{j,s}(1)>0$ for each $j\in [1,n]$ and $s\in [1,r_j-1]$. For each~${i\in [1,n]}$ and $t\in \mathbb{T}_n$, let $F_{i;t}[Z_1,\dots, Z_n]$ be the associated $F$-polynomial of $\mathbf{X}_{i;t}$. We have
\begin{itemize}\itemsep=0pt
 \item[$(1)$] There is a unique monomial $q^{\frac{1}{2}a}\mathbf{Z}(f),a\in\mathbb{Z},$ in $F_{i;t}[Z_1,\dots, Z_n]$ such that it is divisible by all the other monomials in $F_{i;t}[Z_1,\dots, Z_n]$;
 \item[$(2)$] $F_{i;t}[Z_1,\dots, Z_n]$ has constant term $1$.
\end{itemize}
\end{Corollary}

\begin{Remark}
 If $R=I_n$, the assumption in Theorem \ref{t:F-polynomial} is automatically satisfied. Therefore, Theorem \ref{t:F-polynomial}, Corollaries \ref{c:separation-formula} and \ref{c:constant-term-1} reduce to the corresponding results for quantum cluster algebras.
\end{Remark}
\begin{Remark}
 In the setting of cluster algebras and quantum cluster algebras of skew-symmet\-ric type, $F$-polynomials have the positivity property, i.e., each $F$-polynomial has non-negative coefficients (cf.\ \cite{Davison18,GHKK2018}). Thus if one expects that $F$-polynomials for generalized quantum cluster algebras still have the positivity property, then we have to assume that \smash{$h_{i,s}\bigl(q^{\frac{1}{2}}\bigr)$} has non-negative coefficients for each $i$ and $s$, which will imply that $h_{i,s}(1)>0$.
\end{Remark}
\begin{Remark}
 We call the equation \eqref{eq:gupta-formula} Gupta's formula for $F$-polynomials of generalized quantum cluster algebras. When $R=I_n$, it specializes to Gupta's formula
 for quantum cluster algebras. It further specializes to Gupta's formula \cite{gupta2019,LMN2023} for cluster algebras by setting $q^{\frac{1}{2}}=1$. In particular, let $R=I_n$ and set \smash{$q^{\frac{1}{2}}=1$}, we have
 \begin{gather*}
L_{1,i}=\hat{\mathbf{y}}^{\mathbf{c}_{i}^+}\bigl(1+\hat{\mathbf{y}}^{\mathbf{c}_{1}^{+}}\bigr)^{-(d_{(1)}\mathbf{c}_{1},\hat{\mathbf{c}}_{i}^{+})_{D}}\qquad \text{for} \ i\in[1,k],\\
L_{j+1,i}=L_{j,i}(1+L_{j,j+1})^{-(d_{(j+1)}\mathbf{c}_{j+1},\hat{\mathbf{c}}_{i}^{+})_D}\qquad \text{for} \  j\in[1,k-1],\\
L_{1}=1+\hat{\mathbf{y}}^{\mathbf{c}_{1}^{+}},\qquad
L_{l+1}=1+L_{l,l+1}\qquad \text{for} \  l\in[1,k-1].
\end{gather*}
By direct computation,
\begin{align*}
 L_{l+1}={}&1+L_{l-1,l+1}(1+L_{l-1,l})^{-(d_{(l)}\mathbf{c}_{l},\hat{\mathbf{c}}_{l+1}^{+})_D}\\
 ={}&1+L_{l-2,l+1}(1+L_{l-2,l})^{-(d_{(l-1)}\mathbf{c}_{l-1},\hat{\mathbf{c}}_{l+1}^{+})_D}(1+L_{l-1,l})^{-(d_{(l)}\mathbf{c}_{l},\hat{\mathbf{c}}_{l+1}^{+})_D}\\
 & \vdots\\
 ={}&1+L_{1,l+1}\prod\limits^{l}_{p=2}(1+L_{p-1,p})^{-(d_{(p)}\mathbf{c}_{p},\hat{\mathbf{c}}_{l+1}^{+})_D}\\
 ={}&1+\hat{\mathbf{y}}^{\mathbf{c}_{l+1}^{+}}\bigl(1+\hat{\mathbf{y}}^{\mathbf{c}_{1}^{+}}\bigr)^{-(d_{(1)}\mathbf{c}_{1},\hat{\mathbf{c}}_{l+1}^{+})_D}\prod\limits^{l}_{p=2}(1+L_{p-1,p})^{-(d_{(p)}\mathbf{c}_{p},\hat{\mathbf{c}}_{l+1}^{+})_D}\\
={}&1+\hat{\mathbf{y}}^{\mathbf{c}_{l+1}^{+}}\prod\limits^{l}_{p=1}(L_{p})^{-(d_{(p)}\mathbf{c}_{p},\hat{\mathbf{c}}_{l+1}^{+})_D}.
\end{align*}
Hence we obtain Gupta's formula for cluster algebras (cf.\ \cite{LMN2023}):
\smash{$
 F_{i_{k;t_{k}}}=\prod^{k}_{j=1}L_{j}^{-(d_{(j)}\mathbf{c}_{j},\mathbf{g}_{k})_D}$},
where $F_{i_k;t_k}$ is the $i_k$-th $F$-polynomial at vertex $t_k$. We remark the above formula is slightly different from the one in \cite{LMN2023} due to the convention on $\varepsilon_j$.
\end{Remark}

The following example suggests that Theorem \ref{t:F-polynomial}\,(2) may hold without the assumption that~${h_{i,s}(1)>0}$.
\begin{Example}[type $G_{2}$]
Let $\mathcal{A}(3,1)$ be the generalized quantum cluster algebras associated to the initial $(R,\mathbf{h})$-quantum seed $(\mathbf{X},B,\Lambda)$, where
$B=\left[\begin{smallmatrix}
 \hphantom{-}0&1\\
 -1&0
\end{smallmatrix}\right]$ and $\Lambda=\left[\begin{smallmatrix}
 \hphantom{-}0&1\\
 -1&0
\end{smallmatrix}\right]$, $R=\operatorname{diag}\{3,1\}$ and \smash{$\mathbf{h}=\bigl(1,h\bigl(q^{\frac{1}{2}}\bigr),h\bigl(q^{\frac{1}{2}}\bigr),1;1,1\bigr)$},
\smash{$ h\bigl(q^{\frac{1}{2}}\bigr)\in\mathbb{Z}\big[q^{\pm\frac{1}{2}}\big]$}. By assigning $(\mathbf{X}, B,\Lambda)$ to the vertex $t_0\in \mathbb{T}_2$, we obtain an $(R,\mathbf{h})$-quantum seed pattern $t\mapsto \Sigma_t$.

Fix a path \smash{$\xymatrix{t_{0}\ar@{-}[r]^1&t_{1}\ar@{-}[r]^2&t_{2}\ar@{-}[r]^1&t_{3}\ar@{-}[r]^2&t_{4}\ar@{-}[r]^1&t_{5}\ar@{-}[r]^2&t_{6}\ar@{-}[r]^1&t_{7}\ar@{-}[r]^2&t_{8}}$}.
By calculations, we have
\begin{alignat*}{5}
 &\mathbf{g}_{1;t_{1}}=\begin{bmatrix}
 -1\\\hphantom{-}3
 \end{bmatrix}, \qquad&&
 \mathbf{g}_{2;t_{2}}=\begin{bmatrix}
 -1\\\hphantom{-}2
 \end{bmatrix},\qquad&&
 \mathbf{g}_{1;t_{3}}=\begin{bmatrix}
 -2\\\hphantom{-}3
 \end{bmatrix},\qquad&&
  \mathbf{g}_{2;t_{4}}=\begin{bmatrix}
 -1\\\hphantom{-}1
 \end{bmatrix},&\\
 &\mathbf{g}_{1;t_{5}}=\begin{bmatrix}
 -1\\\hphantom{-}0
 \end{bmatrix},\qquad&&
 \mathbf{g}_{2;t_{6}}=\begin{bmatrix}
 \hphantom{-}0\\-1
 \end{bmatrix},\qquad&&
 \mathbf{g}_{1;t_{7}}=\begin{bmatrix}
 1\\0
 \end{bmatrix},\qquad&&
 \mathbf{g}_{2;t_{8}}=\begin{bmatrix}
 0\\1
 \end{bmatrix}.&
 \end{alignat*}
Moreover,
\begin{align*}
 \mathbf{X}_{t_{1}}(e_{1})={}&\mathbf{X}_{t_{0}}(\mathbf{g}_{1;t_1})\bigl(q^{-\frac{3}{2}}\mathbf{\hat{Y}}_{t_0}^{3e_1}+h\bigl(q^{\frac{1}{2}}\bigr)q^{-1}
 \mathbf{\hat{Y}}_{t_0}^{2e_1}+h\bigl(q^{\frac{1}{2}}\bigr)q^{-\frac{1}{2}}\mathbf{\hat{Y}}_{t_0}^{e_1}+1\bigr),\\
 \mathbf{X}_{t_2}(e_2)={}&\mathbf{X}_{t_0}(\mathbf{g}_{2;t_2})\bigl(q^{-\frac{1}{2}}\mathbf{\hat{Y}}_{t_0}^{3e_1+e_2}+q^{-\frac{3}{2}}\mathbf{\hat{Y}}_{t_0}^{3e_1}
 +h\bigl(q^{\frac{1}{2}}\bigr)q^{-1}\mathbf{\hat{Y}}_{t_0}^{2e_1}+h\bigl(q^{\frac{1}{2}}\bigr)q^{-\frac{1}{2}}\mathbf{\hat{Y}}_{t_0}^{e_1}+1\bigr),\\
 \mathbf{X}_{t_3}(e_1)={}&\mathbf{X}_{t_0}(\mathbf{g}_{1;t_3})\bigl(q^{-6}\mathbf{\hat{Y}}_{t_0}^{6e_1}+h\bigl(q^{\frac{1}{2}}\bigr)\bigl(q^{-\frac{11}{2}}+q^{-\frac{9}{2}}\bigr)
 \mathbf{\hat{Y}}_{t_0}^{5e_1}+\big[h\bigl(q^{\frac{1}{2}}\bigr)\bigl(q^{-5}+q^{-3}\bigr)\\
 & +h\bigl(q^{\frac{1}{2}}\bigr)^2q^{-4}\big]\mathbf{\hat{Y}}_{t_0}^{4e_1}+\big[q^{-\frac{9}{2}}+q^{-\frac{3}{2}}+h\bigl(q^{\frac{1}{2}}\bigr)\bigl(q^{-\frac{7}{2}}+q^{-\frac{5}{2}}\bigr)
 \mathbf{\hat{Y}}_{t_0}^{3e_1}\big]\\
 & +[h\bigl(q^{\frac{1}{2}}\bigr)\bigl(q^{-3}+q^{-1}\bigr)+h\bigl(q^{\frac{1}{2}}\bigr)^2q^{-2}\big]\mathbf{\hat{Y}}_{t_0}^{2e_1}+h\bigl(q^{\frac{1}{2}}\bigr)
 \bigl(q^{-\frac{3}{2}}+q^{-\frac{1}{2}}\bigr)\mathbf{\hat{Y}}_{t_0}^{e_1}\\
 & +\bigl(q^{-\frac{11}{2}}+q^{-\frac{9}{2}}+q^{-\frac{7}{2}}\bigr)\mathbf{\hat{Y}}_{t_0}^{6e_1+e_2}+h\bigl(q^{\frac{1}{2}}\bigr)
 \bigl(q^{-\frac{9}{2}}+2q^{-\frac{7}{2}}+q^{-\frac{5}{2}}\bigr)\mathbf{\hat{Y}}_{t_0}^{5e_1+e_2}\\
 & +h\bigl(q^{\frac{1}{2}}\bigr)\bigl(q^{-\frac{7}{2}}+q^{-\frac{5}{2}}+q^{-\frac{3}{2}}+h\bigl(q^{\frac{1}{2}}\bigr)q^{-\frac{5}{2}}\bigr)\mathbf{\hat{Y}}_{t_0}^{4e_1+e_2}\\
 & +\bigl(q^{-\frac{5}{2}}+q^{-\frac{3}{2}}+q^{-\frac{1}{2}}+h\bigl(q^{\frac{1}{2}}\bigr)^2q^{-\frac{3}{2}}\bigr)\mathbf{\hat{Y}}_{t_0}^{3e_1+e_2}
 +h\bigl(q^{\frac{1}{2}}\bigr)q^{-\frac{1}{2}}\mathbf{\hat{Y}}_{t_0}^{2e_1+e_2}\\
 & +\bigl(q^{-4}+q^{-3}+q^{-2}\bigr)\mathbf{\hat{Y}}_{t_0}^{6e_1+2e_2}+h\bigl(q^{\frac{1}{2}}\bigr)\bigl(q^{-\frac{5}{2}}+q^{-\frac{3}{2}}\bigr)\mathbf{\hat{Y}}_{t_0}^{5e_1+2e_2}\\
 & +h\bigl(q^{\frac{1}{2}}\bigr)q^{-1}\mathbf{\hat{Y}}_{t_0}^{4e_1+2e_2}+q^{-\frac{3}{2}}\mathbf{\hat{Y}}_{t_0}^{6e_1+3e_2}+1
 \bigr),\\
 \mathbf{X}_{t_4}(e_2)={}&\mathbf{X}_{t_0}(\mathbf{g}_{2;t_4})\bigl(q^{-\frac{3}{2}}\mathbf{\hat{Y}}_{t_0}^{3e_1}+h\bigl(q^{\frac{1}{2}}\bigr)q^{-1}\mathbf{\hat{Y}}_{t_0}^{2e_1}
 +h\bigl(q^{\frac{1}{2}}\bigr)q^{-\frac{1}{2}}\mathbf{\hat{Y}}_{t_0}^{e_1}\\
 & +\bigl(q^{-\frac{3}{2}}+q^{-\frac{1}{2}}\bigr)\mathbf{\hat{Y}}_{t_0}^{3e_1+e_2}+h\bigl(q^{\frac{1}{2}}\bigr)q^{-\frac{1}{2}}\mathbf{\hat{Y}}_{t_0}^{2e_1+e_2}
 +q^{-\frac{1}{2}}\mathbf{\hat{Y}}_{t_0}^{3e_1+2e_2}+1
 \bigr),\\
 \mathbf{X}_{t_5}(e_1)={}&\mathbf{X}_{t_0}(\mathbf{g}_{1;t_5})\bigl(q^{-\frac{3}{2}}\mathbf{\hat{Y}}_{t_0}^{3e_1}+h\bigl(q^{\frac{1}{2}}\bigr)q^{-1}\mathbf{\hat{Y}}_{t_0}^{2e_1}+h\bigl(q^{\frac{1}{2}}\bigr)q^{-\frac{1}{2}}\mathbf{\hat{Y}}_{t_0}^{e_1}\\
 & +\bigl(q^{-\frac{5}{2}}+q^{-\frac{3}{2}}+q^{-\frac{1}{2}}\bigr)\mathbf{\hat{Y}}_{t_0}^{3e_1+e_2}+h\bigl(q^{\frac{1}{2}}\bigr)\bigl(q^{-\frac{3}{2}}+q^{-\frac{1}{2}}\bigr)
 \mathbf{\hat{Y}}_{t_0}^{2e_1+e_2}+h\bigl(q^{\frac{1}{2}}\bigr)q^{-\frac{1}{2}}\mathbf{\hat{Y}}_{t_0}^{e_1+e_2}\\
 & +\bigl(q^{-\frac{5}{2}}+q^{-\frac{3}{2}}+q^{-\frac{1}{2}}\bigr)\mathbf{\hat{Y}}_{t_0}^{3e_1+2e_2}+h\bigl(q^{\frac{1}{2}}\bigr)q^{-1}\mathbf{\hat{Y}}_{t_0}^{2e_1+2e_2}
 +q^{-\frac{3}{2}}\mathbf{\hat{Y}}_{t_0}^{3e_1+3e_2}+1
 \bigr),\\
 \mathbf{X}_{t_6}(e_2)={}&\mathbf{X}_{t_0}(\mathbf{g}_{2;t_6})\bigl(q^{-\frac{1}{2}}\mathbf{\hat{Y}}_{t_0}^{e_2}+1\bigr),\qquad
 \mathbf{X}_{t_7}(e_1)= \mathbf{X}_{t_0}(\mathbf{g}_{1;t_7})=\mathbf{X}_{t_0}(e_1),\\
 \mathbf{X}_{t_8}(e_2)={}&\mathbf{X}_{t_0}(\mathbf{g}_{2;t_8})=\mathbf{X}_{t_0}(e_{2}),
\end{align*}
where $e_1$, $e_2$ is the standard $\mathbb{Z}$-basis of $\mathbb{Z}^2$.
In particular, each $F_{i_j;t_j}$ is a polynomial in $\mathbf{\hat{Y}}_{1}$, $\mathbf{\hat{Y}}_2$ for any choice of \smash{$h\bigl(q^{\frac{1}{2}}\bigr)$}. Moreover, it has non-negative coefficients when regarding \smash{$h\bigl(q^{\frac{1}{2}}\bigr)$} as a~variable. This phenomenon will be studied in a forthcoming paper.
\end{Example}

\section{Gupta's formula for generalized cluster algebras} \label{s:gupta-formula-gca}
As we have seen in the proof of Theorem \ref{t:F-polynomial} that $F$-polynomials of a generalized quantum cluster algebra do not specialize to $F$-polynomials of a generalized cluster algebra with principal coefficients, hence Gupta's formula \eqref{eq:gupta-formula} does not reduce to a formula of $F$-polynomials for a~generalized cluster algebra.
However, we can prove the Gupta's formula for $F$-polynomials of generalized cluster algebras following Sections \ref{S5} and \ref{S6} (cf.\ \cite{LMN2023}).

Keep notation as in Section \ref{ss:principal-gca}. Recall that $(\mathbf{x},\mathbf{y},B)$ is a labeled seed in $\mathcal{F}$ with coefficients in $\mathbb{P}=\operatorname{Trop}(\mathbf{y},\mathbf{z})$.
By assigning $(\mathbf{x},\mathbf{y},B)$ to the vertex $t_0\in \mathbb{T}_n$, we obtain an $(\mathbf{r},\mathbf{z})$-seed pattern~${t\mapsto \Sigma_t}$.

\subsection{Fock--Goncharov decomposition}\label{S3}
In this section, we establish Fock--Goncharov decomposition for generalized cluster algebras with principal coefficients, which generalizes the corresponding construction for cluster algebras~\mbox{\cite{fock2009,nakanishi2023}}.

Let \smash{$\xymatrix{t\ar@{-}[r]^k&t'}$} in $\mathbb{T}_n$ and $k\in [1,n]$, denote by $\varepsilon_{k;t}$ the common sign of components of the $c$-vector $\mathbf{c}_{k;t}$. As before, we also set $\mathbf{c}_{k;t}^+:=\varepsilon_{k;t}\mathbf{c}_{k;t}$ and $\mathbf{\hat{c}}_{k;t}=B\mathbf{c}_{k;t}$. By the sign-coherence of~$\mathbf{c}_{k;t}$, the exchange formula \eqref{E4} at $x_{k;t}$ can be simplified as
\[
x_{k;t'}=x_{k;t}^{-1}\Biggl(\prod\limits_{j=1}^{n}x_{j;t}^{[-\varepsilon_{k;t} b_{jk;t}]_+}\Biggr)^{r_k}\sum\limits_{s=0}^{r_k}z_{k,s}\hat{y}_{k;t}^{\varepsilon_{k;t} s}.
\]
Let $\mathbb{Q}\mathbb{P}(\mathbf{x}_t)$ be the field of rational functions in $\mathbf{x}_t$ with coefficients in $\mathbb{Q}\mathbb{P}$. The mutation $\mu_{k;t}$ yields an isomorphism $\mu_{k;t}\colon\mathbb{Q}\mathbb{P}(\mathbf{x}_{t'})\to \mathbb{Q}\mathbb{P}(\mathbf{x}_t)$ which maps $x_{i;t'}$ to $x_{i;t}$ for $i\neq k$ and maps~$x_{k;t'}$~to
\[
x_{k;t}^{-1}\Biggl(\prod\limits_{j=1}^{n}x_{j;t}^{[-\varepsilon_{k;t} b_{jk;t}]_+}\Biggr)^{r_k}\left(\sum\limits_{s=0}^{r_k}z_{k,s}\hat{y}_{k;t}^{\varepsilon_{k;t} s}\right).
\]

 We introduce two homomorphisms of fields as follows:
\begin{align*}
 \tau_{k;t}\colon\ \mathbb{QP}(\mathbf{x}_t')\longrightarrow\mathbb{QP}(\mathbf{x}_t),\qquad
 x_{i;t'}\mapsto
 \begin{cases}
 x_{i;t} & \text{if $i\neq k$,}\\
 \displaystyle x_{k;t}^{-1}\biggl(\prod\limits_{j=1}^{n}x_{j;t}^{[-\varepsilon_{k;t}b_{jk;t}r_k]_+}\biggr) & \text{if $i=k$,}
 \end{cases}
 \end{align*} and
 \begin{align*}
 \rho_{k;t}\colon\ \mathbb{QP}(\mathbf{x}_t)\longrightarrow\mathbb{QP}(\mathbf{x}_t),\qquad
 x_{i;t}\mapsto x_{i;t}\left(\sum\limits_{s=0}^{r_k}z_{k,s}\hat{y}^{\varepsilon_{k;t}s}_{k;t}\right)^{-\delta_{i,k}}.
 \end{align*}
It is clear that $\mu_{k;t}=\rho_{k;t}\circ\tau_{k;t}$ and $\mu_{k;t}\circ\mu_{k;t'}=\id$. Moreover, we have the following by direct computation.
\begin{Lemma}\label{P8}
 The following relation holds:
$
 \tau_{k;t}\circ\tau_{k;t'}=\id$.
\end{Lemma}

We further introduce an automorphism $q_{k;t}$ of $\mathbb{QP}(\mathbf{x})$ as follows
\begin{align*}
 q_{k;t}\colon\
 \mathbb{QP}(\mathbf{x})\longrightarrow\mathbb{QP}(\mathbf{x}),\qquad
 \mathbf{x}^{\mathbf{m}}\longmapsto\mathbf{x}^{\mathbf{m}}\left(\sum\limits_{s=0}^{r_k}z_{k,s}\bigl(\hat{\mathbf{y}}^{\mathbf{c}_{k;t}^+}\bigr)^s\right)^{-(\mathbf{m},d_{k}\mathbf{c}_{k;t})_{D_0R}},\qquad \forall \mathbf{m} \in\mathbb{Z}^n.
\end{align*}
By computation, we have
\begin{align}
 q_{k;t}(\hat{\mathbf{y}}^{\mathbf{n}})=\hat{\mathbf{y}}^{\mathbf{n}}
 \left(\sum\limits_{s=0}^{r_k}z_{k,s}\bigl(\hat{\mathbf{y}}^{\mathbf{c}_{k;t}^+}\bigr)^s\right)^{(\mathbf{n},d_{k}\mathbf{\hat{c}}_{k;t})_{D_0R}},\qquad \forall \mathbf{n}\in\mathbb{Z}^n.\label{E22}
\end{align}
\begin{Lemma}\label{P10}
The following formulas hold:
 \begin{gather}
  q_{k;t}\bigl(\mathbf{x}^{\mathbf{g}_{i;t}}\bigr)=\mathbf{x}^{\mathbf{g}_{i;t}}\left(\sum\limits_{s=0}^{r_k}z_{k,s}\bigl(\hat{\mathbf{y}}^{\mathbf{c}_{k;t}^+}\bigr)^s\right)^{-\delta_{i;k}},\label{E23}\\
  q_{k;t'}\circ q_{k;t}=\id. \label{E25}
 \end{gather}
\end{Lemma}
\begin{proof}
Formula \eqref{E23} is a direct consequence of \eqref{E12}.

 For any $\mathbf{m}\in\mathbb{Z}^n$,
 \begin{align*}
 q_{k;t'}\circ q_{k;t}(\mathbf{x}^{\mathbf{m}})={}&q_{k;t'}
 \left(\mathbf{x}^{\mathbf{m}}\left(\sum\limits_{s=0}^{r_k}z_{k,s}\bigl(\hat{\mathbf{y}}^{\mathbf{c}_{k;t}^+}\bigr)^s\right)^{-(\mathbf{m},d_{k}\mathbf{c}_{k;t})_{D_0R}}\right)\\
 ={}&\mathbf{x}^{\mathbf{m}}\left(\sum\limits_{s=0}^{r_k}z_{k,s}\bigl(\hat{\mathbf{y}}^{\mathbf{c}_{k;t'}^+}\bigr)^s\right)^{-(\mathbf{m},d_{k}\mathbf{c}_{k;t'})_{D_0R}}\\
 &\times\left(\sum\limits_{s=0}^{r_k}z_{k,s}\left(\hat{\mathbf{y}}^{\mathbf{c}_{k;t}^+}
 \left(\sum\limits_{s=0}^{r_k}z_{k,s}\bigl(\hat{\mathbf{y}}^{\mathbf{c}_{k;t'}^+}\bigr)^s\right)^{(\mathbf{c}_{k;t}^+,d_k\mathbf{\hat{c}}_{k;t'})_{D_0R}}\right)^s\right)^{-(\mathbf{m},d_k\mathbf{c}_{k;t})_{D_0R}}.
 \end{align*}
 Note that $D_0RB$ is skew-symmetric, \smash{$\bigl(\mathbf{c}_{k;t}^+,d_k\mathbf{\hat{c}}_{k;t'}\bigr)_{D_0R}=-d_k\varepsilon_{k;t}(\mathbf{c}_{k;t})^{\mathsf T}D_0RB\mathbf{c}_{k;t}=0$}. On the other hand, $\mathbf{c}_{k;t'}=-\mathbf{c}_{k;t}$. Putting all of these together, we obtain that
 $
 q_{k;t'}\circ q_{k;t}\bigl(\mathbf{x}^\mathbf{m}\bigr)=\mathbf{x}^\mathbf{m}$.
 This completes the proof of \eqref{E25}.
\end{proof}

From now on, we fix a path
\smash{$\xymatrix{t_{0}\ar@{-}[r]^{i_1}&t_{1}\ar@{-}[r]^{i_2}&t_{2}\ar@{-}[r]^{i_3}&\cdots\ar@{-}[r]^{i_k}&t_{k}}$} in $\mathbb{T}_n$. We define
\begin{align*}
 &\mu_{t_k}^{t_0}:=\mu_{i_1;t_0}\circ \mu_{i_2;t_1}\circ\cdots\circ \mu_{i_k;t_{k-1}}\colon\ \mathbb{QP}(\mathbf{x}_{t_k})\to \mathbb{QP}(\mathbf{x}_{t_0}),\\
 &\tau_{t_k}^{t_0}:=\tau_{i_1;t_0}\circ\tau_{i_2;t_1}\circ\cdots\circ \tau_{i_k;t_{k-1}}\colon\  \mathbb{QP}(\mathbf{x}_{t_k})\to \mathbb{QP}(\mathbf{x}_{t_0}),\\
 &q_{t_k}^{t_0}:=q_{i_1;t_0}\circ q_{i_2;t_1}\circ\cdots\circ q_{i_k;t_{k-1}}\colon\   \mathbb{QP}(\mathbf{x}_{t_0})\to \mathbb{QP}(\mathbf{x}_{t_0}).
\end{align*}

\begin{Lemma}\label{L9}
For any $i\in [1,n]$, the following formulas hold:
\begin{align}
& \mu_{t_k}^{t_0}(x_{i;t_k})=\mathbf{x}^{\mathbf{g}_{i;t_k}}F_{i;t_k}(\mathbf{\hat{y}},\mathbf{z}),\label{E19}\\
 &\tau_{t_k}^{t_0}(x_{i;t_k})=\mathbf{x}^{\mathbf{g}_{i;t_k}},\label{E20}\\
 &\tau_{t_k}^{t_0}(\hat{y}_{i;t_k})=\mathbf{\hat{y}}^{\mathbf{c}_{i;t_k}}\label{E21}.
\end{align}
\end{Lemma}
\begin{proof}
Formula \eqref{E19} is equivalent to the separation formula of generalized cluster alge\-bras~\cite{Nakanishi14}. Formula \eqref{E20} follows form \eqref{E15} by induction on $k$, and \eqref{E21} follows from \eqref{E20} and~\eqref{E14}.
\end{proof}

\begin{Lemma}
 For any $j\in [1,k]$ and $l\in [1,n]$, we have
 \begin{align}
 \tau_{t_j}^{t_0}\circ\rho_{l;t_j}&=q_{l;t_j}\circ\tau_{t_j}^{t_0}.\label{E24}
\end{align}
\end{Lemma}
\begin{proof}
For any $i\in [1,n]$,
\begin{align*}
 \tau_{t_j}^{t_0}\circ\rho_{l;t_j}(x_{i;t_j})
 &=\tau_{t_j}^{t_0}\left(x_{i;t_j}\left(\sum\limits_{s=0}^{r_{l}}z_{l,s}\bigl(\hat{y}^{\varepsilon_{l;t_j}}_{l;t_j}\bigr)^s\right)^{-\delta_{i,l}}\right)
 =\mathbf{x}^{\mathbf{g}_{i;t_j}}\left(\sum\limits_{s=0}^{r_{l}}z_{l,s}\bigl(\mathbf{\hat{y}}^{\varepsilon_{l;t_j}\mathbf{c}_{l;t_j}}\bigr)^s\right)^{-\delta_{i,l}}\\
 &=\mathbf{x}^{\mathbf{g}_{i;t_j}}\left(\sum\limits_{s=0}^{r_{l}}z_{l,s}\bigl(\mathbf{\hat{y}}^{\mathbf{c}_{l;t_j}^+}\bigr)^s\right)^{-\delta_{i,l}}
 =q_{l;t_j}\circ\tau_{t_j}^{t_0}(x_{i;t_j}),
 \end{align*}
 where the last equality follows from \eqref{E23} and \eqref{E20}.
\end{proof}

By applying \eqref{E24}, we have
\begin{Proposition}\label{P11}
 The following decomposition holds:
$
\mu_{t_k}^{t_0}=q_{t_k}^{t_0}\circ\tau_{t_k}^{t_0}\colon\  \mathbb{QP}(\mathbf{x}_{t_k})\longrightarrow\mathbb{QP}(\mathbf{x}_{t_0})$.
\end{Proposition}

The following is a direct consequence of Proposition \ref{P11} and Lemma \ref{L9}.
\begin{Proposition}\label{p:FG-decomp}
 For any $i\in [1,n]$, the following formula holds:
$
 q_{t_k}^{t_0}\bigl(\mathbf{x}^{\mathbf{g}_{i;t_k}}\bigr)=\mathbf{x}^{\mathbf{g}_{i;t_k}}F_{i;t_k}(\hat{\mathbf{y}},\mathbf{z})$.
\end{Proposition}

\subsection{Gupta's formula}\label{S5}
Throughout this section, we fix a path
\smash{$\xymatrix{t_{0}\ar@{-}[r]^{i_1}&t_{1}\ar@{-}[r]^{i_2}&t_{2}\ar@{-}[r]^{i_3}&\cdots\ar@{-}[r]^{i_k}&t_{k}}$}
 in $\mathbb{T}_n$. For the simplicity of notation, for any $j\in [1,k]$, we denote by \smash{$d_{(j)}:=d_{i_j}$}, \smash{$r_{(j)}:=r_{i_j}$}, \smash{$\mathbf{c}_j:=\mathbf{c}_{i_j;t_{j-1}}$}, \smash{$\mathbf{c}_j^+\!\!:=\varepsilon_{i_j;t_{j-1}}\mathbf{c}_{j}$}, \smash{$\hat{\mathbf{c}}_j^+\!:=\!B\mathbf{c}_j^+$}, \smash{$\mathbf{g}_j\!:=\mathbf{g}_{i_j;t_j}$}, where $\varepsilon_{i_j;t_{j-1}}$ is the common sign of components of~$\mathbf{c}_{i_j;t_{j-1}}$.
We also introduce certain elements $L_1,\dots, L_k$ of $\mathcal{F}$ along the path, where
\[
L_1=\sum_{s=1}^{r_{(1)}}\!z_{i_1,s}\bigl(\hat{\mathbf{y}}^{\mathbf{c}_{1}^+}\bigr)^s, \qquad
L_l=\sum\limits_{s=0}^{r_{(l)}}z_{i_l,s}\left(\hat{\mathbf{y}}^{\mathbf{c}_{l}^+} \prod\limits_{j=1}^{l-1}L_j^{-(\mathbf{\hat{c}}_{l}^+,d_{(j)}\mathbf{c}_{j})_{D_0R}}\right)^s\qquad \text{for}\  2\leq l\leq k.
\]

The following is the Gupta's formula of $F$-polynomials for generalized cluster algebras.
\begin{Theorem}\label{T5}
The following formula holds:
 \begin{align}
 F_{i_k;t_k}(\hat{\mathbf{y}},\mathbf{z})=\prod\limits_{j=1}^{k}L_{j}^{-(\mathbf{g}_{k},d_{(j)}\mathbf{c}_{j})_{D_0R}}.\label{E34}
 \end{align}
\end{Theorem}

Before giving the proof of Theorem \ref{T5}, we prepare the following lemma.
\begin{Lemma}\label{L3}
For $1\leq i\leq m\leq k$, the following formula holds:
 \begin{align}\label{eq:q_t-formula}
 q_{t_i}^{t_0}\bigl(\hat{\mathbf{y}}^{\mathbf{c}_{m}^+}\bigr)=\hat{\mathbf{y}}^{\mathbf{c}_{m}^+}
 \Biggl(\prod\limits_{j=1}^{i}L_j^{-(\mathbf{\hat{c}}_{m}^+,d_{(j)}\mathbf{c}_{j})_{D_0R}}\Biggr).
 \end{align}

\end{Lemma}
\begin{proof}
 We prove this formula by induction on $i$.
 For $i=1$, by \eqref{E22}, we get
\begin{align*}
 q_{t_{1}}^{t_0}\bigl(\hat{\mathbf{y}}^{\mathbf{c}_{m}^+}\bigr)&=\hat{\mathbf{y}}^{\mathbf{c}_{m}^+}
 \left(\sum\limits_{s=0}^{r_{(1)}}z_{i_1,s}\bigl(\hat{\mathbf{y}}^{\mathbf{c}_{1}^+}\bigr)^s\right)^{(\mathbf{c}_{m}^+,d_{(1)}\mathbf{\hat{c}}_{i_{1};t_{0}})_{D_0R}}=\hat{\mathbf{y}}^{\mathbf{c}_{m}^+}(L_1)^{-(\mathbf{\hat{c}}_{m}^+,d_{(1)}\mathbf{c}_{1})_{D_0R}}.
\end{align*}
Now suppose that
\begin{align*}
 q_{t_{l-1}}^{t_0}\bigl(\hat{\mathbf{y}}^{\mathbf{c}_{m}^+}\bigr)=\hat{\mathbf{y}}^{\mathbf{c}_{m}^+}
 \left(\prod\limits_{j=1}^{l-1}L_j^{-(\mathbf{\hat{c}}_{m}^+,d_{(j)}\mathbf{c}_{j})_{D_0R}}\right).
\end{align*}
Then,
\begin{align*}
 q_{t_l}^{t_0}\bigl(\hat{\mathbf{y}}^{\mathbf{c}_{m}^+}\bigr)={}&q_{t_{l-1}}^{t_0}\bigl(q_{i_l;t_{l-1}}\bigl(\hat{\mathbf{y}}^{\mathbf{c}_{m}^+}\bigr)\bigr)
 =q_{t_{l-1}}^{t_0}\Biggl(\hat{\mathbf{y}}^{\mathbf{c}_{m}^+}
 \left(\sum\limits_{s=0}^{r_{(l)}}z_{i_l,s}(\hat{\mathbf{y}}^{\mathbf{c}_{l}^+})^s\right)^{-(\mathbf{\hat{c}}_{m}^+,d_{(l)}\mathbf{c}_{l})_{D_0R}}\Biggr)\\
 ={}&\hat{\mathbf{y}}^{\mathbf{c}_{m}^+}\left(\prod\limits_{j=1}^{l-1}L_j^{-(\mathbf{\hat{c}}_{m}^+,d_{(j)}\mathbf{c}_{j})_{D_0R}}\right)\\
&\times\left(\sum\limits_{s=0}^{r_{(l)}}z_{i_l,s}\hat{\mathbf{y}}^{s \mathbf{c}_{l}^+}
\left(\prod\limits_{j=1}^{l-1}L_j^{-s(\mathbf{\hat{c}}_{l}^+,d_{(j)}\mathbf{c}_{j})_{D_0R}}\right)\right)^{-(\mathbf{\hat{c}}_{m}^+,d_{(l)}\mathbf{c}_{l})_{D_0R}}\\
={}&\hat{\mathbf{y}}^{\mathbf{c}_{m}^+}\left(\prod\limits_{j=1}^{l-1}L_j^{-(\mathbf{\hat{c}}_{m}^+,d_{(j)}\mathbf{c}_{j})_{D_0R}}\right)
(L_l)^{-(\mathbf{\hat{c}}_{m}^+,d_{(l)}\mathbf{}c_{l})_{D_0R}}\\
={}&\hat{\mathbf{y}}^{\mathbf{c}_{m}^+}\left(\prod\limits_{j=1}^{l}L_j^{-(\mathbf{\hat{c}}_{m}^+,d_{(j)}\mathbf{c}_{j})_{D_0R}}\right).\tag*{\qed}
\end{align*} \renewcommand{\qed}{}
\end{proof}

{\bf Proof of Theorem \ref{T5}.}
For $2\leq l\leq k$, we have
\begin{align*}
 \frac{q_{t_{l}}^{t_0}(\mathbf{x}^{\mathbf{g}_k})}{q_{t_{l-1}}^{t_0}(\mathbf{x}^{\mathbf{g}_k})}
 &=q_{t_{l-1}}^{t_0}\left(\frac{q_{i_l;t_{l-1}}(\mathbf{x}^{\mathbf{g}_k})}{\mathbf{x}^{\mathbf{g}_k}}\right)\\
 &=q_{t_{l-1}}^{t_0}\Biggl(\left(\sum\limits_{s=0}^{r_{(l)}}z_{i_l,s}\bigl(\hat{\mathbf{y}}^{\mathbf{c}_{l}^+}\bigr)^s\right)^{-(\mathbf{g}_{k},d_{(l)}\mathbf{c}_{l})_{D_0R}}\Biggr)\\
 &=\left(\sum\limits_{s=0}^{r_{(l)}}z_{i_l,s}\bigl(q_{t_{l-1}}^{t_0}\bigl(\hat{\mathbf{y}}^{\mathbf{c}_{l}^+}\bigr)\bigr)^s\right)^{-(\mathbf{g}_{k},d_{(l)}\mathbf{c}_{l})_{D_0R}}\\
 &=\left(\sum\limits_{s=0}^{r_{(l)}}z_{i_l,s}\left(\hat{\mathbf{y}}^{\mathbf{c}_{l}^+}\prod\limits_{j=1}^{l-1}L_j^{-(\mathbf{\hat{c}}_{l}^+,d_{(j)}\mathbf{c}_{j})_{D_0R}}\right)^s\right)^{-(\mathbf{g}_{k},d_{(l)}\mathbf{c}_{l})_{D_0R}} \qquad (\text{by \eqref{eq:q_t-formula}})\\
 &=(L_l)^{-(\mathbf{g}_{k},d_{(l)}\mathbf{c}_{l})_{D_0R}}.
\end{align*}
On the other hand, \[\frac{q_{i_1;t_0}(\mathbf{x}^{\mathbf{g}_{k}})}{\mathbf{x}^{\mathbf{g}_{k}}}=\left(\sum\limits_{s=0}^{r_{(1)}}z_{i_1,s}\bigl(\hat{\mathbf{y}}^{\mathbf{c}_{1}^+}\bigr)^s\right)^{-(\mathbf{g}_{k},d_{(1)}\mathbf{c}_{1})_{D_0R}}=L_1^{-(\mathbf{g}_{k},d_{(1)}\mathbf{c}_{1})_{D_0R}}.\]
According to Proposition \ref{p:FG-decomp}, we have
\begin{align*}
 F_{i_k;t_k}(\hat{\mathbf{y}},\mathbf{z})&=\frac{q_{t_k}^{t_0}(\mathbf{x}^{\mathbf{g}_{k}})}{\mathbf{x}^{\mathbf{g}_k}}=\frac{q_{t_k}^{t_0}(\mathbf{x}^{\mathbf{g}_{k}})}{q_{t_{k-1}}^{t_0}(\mathbf{x}^{\mathbf{g}_{k}})}\circ\frac{q_{t_{k-1}}^{t_0}(\mathbf{x}^{\mathbf{g}_{k}})}{q_{t_{k-2}}^{t_0}(\mathbf{x}^{\mathbf{g}_{k}})}\circ\cdots\circ\frac{q_{t_1}^{t_0}(\mathbf{x}^{\mathbf{g}_{k}})}{\mathbf{x}^{\mathbf{g}_{k}}}\\
 &=(L_k)^{-(\mathbf{g}_{k},d_{(k)}\mathbf{c}_{k})_{D_0R}}(L_{k-1})^{-(\mathbf{g}_{k},d_{(k-1)}\mathbf{c}_{k-1})_{D_0R}}\cdots(L_1)^{-(\mathbf{g}_{k},d_{(1)}\mathbf{c}_{1})_{D_0R}}\\
 &=\prod\limits_{j=1}^{k}L_j^{-(\mathbf{g}_{k},d_{(j)}\mathbf{c}_{j})_{D_0R}}.
\end{align*}
This completes the proof. \hfill \qedsymbol

\begin{Example}
\label{EX1}
Let us explain Theorem \ref{T5} by the following simplest non-trivial example. Let~${\mathbf{r}=(2,1)}$, $\mathbf{z}=(1,z,1;1,1)$ and $B=\left[\begin{smallmatrix}
 \hphantom{-}0 & 1\\
 -1 & 0
 \end{smallmatrix}\right]$. We consider the following path in $\mathbb{T}_2$:
\smash{$
 \xymatrix{t_0\ar@{-}[r]^1&t_1\ar@{-}[r]^2&t_2\ar@{-}[r]^1&t_3\ar@{-}[r]^2&t_4}$}.
 By assigning the labeled seed $\Sigma=(\mathbf{x},\mathbf{y},B)$ to the vertex~$t_0$, we obtain an $(\mathbf{r},\mathbf{z})$-seed pattern~${t\mapsto \Sigma_t}$.
By definition, $\hat{y}_1=y_1x_2^{-1}$, $\hat{y}_2=y_2x_1$, $RB_{t_0}=\left[\begin{smallmatrix}
 \hphantom{-}0 & 2 \\
 -1 & 0
 \end{smallmatrix}\right]$, $d_{(1)}=d_{(2)}=d_{(3)}=d_{(4)}=\frac{1}{2}$, and the $c$-,$\hat{c}$- and $g$-vectors involved are as following:
 \begin{alignat*}{5}
 &\mathbf{c}_{1}=\begin{bmatrix}
 1\\0
 \end{bmatrix}, \qquad&&
 \mathbf{c}_{2}=\begin{bmatrix}
 2\\1
 \end{bmatrix},\qquad&&
 \mathbf{c}_{3}=\begin{bmatrix}
 1\\1
 \end{bmatrix},\qquad&&
 \mathbf{c}_{4}=\begin{bmatrix}
 0\\1
 \end{bmatrix},&\\
 &\mathbf{\hat{c}}_{1}=\begin{bmatrix}
 \hphantom{-}0\\-1
 \end{bmatrix},\qquad&&
 \mathbf{\hat{c}}_{2}=\begin{bmatrix}
 \hphantom{-}1\\-2
 \end{bmatrix},\qquad&&
 \mathbf{\hat{c}}_{3}=\begin{bmatrix}
 \hphantom{-}1\\-1
 \end{bmatrix},\qquad&&
 \mathbf{\hat{c}}_{4}=\begin{bmatrix}
 1\\0
 \end{bmatrix},&\\
 &\mathbf{g}_{1}=\begin{bmatrix}
 -1\\\hphantom{-}2
 \end{bmatrix},\qquad&&
 \mathbf{g}_{2}=\begin{bmatrix}
 -1\\\hphantom{-}1
 \end{bmatrix},\qquad&&
 \mathbf{g}_{3}=\begin{bmatrix}
 -1\\\hphantom{-}0
 \end{bmatrix},\qquad&&
 \mathbf{g}_{4}=\begin{bmatrix}
 \hphantom{-}0\\-1
 \end{bmatrix}.&
 \end{alignat*}
The relevant inner products\footnote{Here we omit the subscript $D_0R$.} are given by
\begin{alignat*}{4}
&\bigl(\mathbf{\hat{c}}_2^+,d_{(1)}\mathbf{c}_{1}\bigr)=1,\qquad &&\bigl(\mathbf{\hat{c}}_3^+,d_{(1)}\mathbf{c}_{1}\bigr)=1,\qquad && \bigl(\mathbf{\hat{c}}_3^+,d_{(2)}\mathbf{c}_{2}\bigr)=1,&\\
 &\bigl(\mathbf{\hat{c}}_4^+,d_{(1)}\mathbf{c}_{1}\bigr)=1,\qquad&&
 \bigl(\mathbf{\hat{c}}_4^+,d_{(2)}\mathbf{c}_{2}\bigr)=2,\qquad&&
\bigl(\mathbf{\hat{c}}_4^+,d_{(3)}\mathbf{c}_{3}\bigr)=1,&\\
 &(\mathbf{g}_{2},d_{(1)}\mathbf{c}_{1})=-1,\qquad&&
 (\mathbf{g}_{3},d_{(1)}\mathbf{c}_{1})=-1,\qquad&&
 (\mathbf{g}_{3},d_{(2)}\mathbf{c}_{2})=-2,&\\
 &(\mathbf{g}_{4},d_{(1)}\mathbf{c}_{1})=0,\qquad&&
 (\mathbf{g}_{4},d_{(2)}\mathbf{c}_{2})=-1,\qquad&&
 (\mathbf{g}_{4},d_{(3)}\mathbf{c}_{3})=-1.&
 \end{alignat*}
 Therefore,
 \begin{align*}
 L_1&=1+z\hat{y}_1+\hat{y}_1^2,\qquad
 L_2=1+\hat{y}_1^2\hat{y}_2L_1^{-1},\\
 L_3&=1+z\bigl(\hat{y}_1\hat{y}_2L_1^{-1}L_2^{-1}\bigr)+\bigl(\hat{y}_1\hat{y}_2L_1^{-1}L_2^{-1}\bigr)^2,\qquad
 L_4=1+\hat{y}_2L_1^{-1}L_2^{-2}L_{3}^{-1}.
 \end{align*}
 Applying Theorem \ref{T5}, we get
 \begin{align*}
 F_{1;t_1}&=L_1
 =1+z\hat{y}_1+\hat{y}_1^2,\qquad
 F_{2;t_2}=L_1L_2=1+z\hat{y}_1+\hat{y}_1^2+\hat{y}_1^2\hat{y}_2,\\
 F_{1;t_3}&=L_1L_2^2L_3=1+z\hat{y}_1+\hat{y}_1^2+z\hat{y}_1\hat{y}_2+2\hat{y}_1^2\hat{y}_2+\hat{y}_1^2\hat{y}_2^2,\\
 F_{2;t_4}&=L_2L_3L_4
 =1+\hat{y}_2.
\end{align*}
\end{Example}

\subsection{Expansion of Gupta's formula}\label{S6}
Following \cite{LMN2023}, we expand Gupta's formula \eqref{E34} into sum in this section. Let $h\in \Z$ and $n_0,\dots, n_l\in \N$, we denote
\[\begin{Bmatrix}
 h\\
 n_0,n_1,\dots, n_l
\end{Bmatrix}:=\binom{h}{n_0}\binom{n_0}{n_1,\dots, n_l},
\]
where $\binom{n_0}{n_1,\dots, n_l}$ is the multinomial coefficient. In particular, if $n_0\neq n_1+\cdots +n_l$, then \smash{$\left\{\begin{smallmatrix}
 h\\
 n_0,n_1,\dots, n_l
\end{smallmatrix}\right\}=0$}.
It easy to see that if $h>0$, then \[\begin{Bmatrix}
 h\\
 n_0,n_1,\dots, n_l
\end{Bmatrix}=\binom{h}{h-n_0,n_1,\dots, n_l}.\]
\begin{Lemma}\label{expansion of polynomial}
 Let $1+\sum_{i=1}^{l}a_iz^i$ be a polynomial in $z$ with coefficients in a field $\mathbb{F}$ and $h\in \Z$, we have the following expansion formula:
 \[
 \left(1+\sum\limits_{i=1}^{l}a_iz^i\right)^h=\sum\begin{Bmatrix}
 h\\
 n_0,n_1,\dots, n_l
\end{Bmatrix}a_1^{n_1}\cdots a_l^{n_l}z^{\sum\limits_{i=1}^{l} in_i},\]
where the sum takes over $n_0,\dots, n_m \geq 0$ and $n_0=n_1+\cdots+n_l$.
\end{Lemma}
\begin{proof}
 Denote by $w$ the polynomial $\sum_{i=1}^{l}a_iz^i$. For any integer $h$, we have \[(1+w)^{h}=\sum\limits_{n_0\ge 0}\binom{h}{n_0}w^{n_0}.\]
 Then replacing $w$ with $\sum_{i=1}^{l}a_iz^i$, we obtain
 \begin{align*}
 \left(1+\sum\limits_{i=1}^{l}a_iz^i\right)^h=\sum\binom{h}{n_0}\binom{n_0}{n_1,\dots,n_l}a_1^{n_1}\cdots a_l^{n_l}z^{\sum\limits_{i=1}^{l} in_i}.\tag*{\qed}
 \end{align*}\renewcommand{\qed}{}
\end{proof}

\begin{Lemma}\label{L4}
 Given integers $h_1,\dots,h_l$ and the same setup as Theorem {\rm\ref{T5}},
\begin{align*} \prod\limits_{j=1}^{k}L_j^{h_j}=\sum\prod\limits_{j=1}^{k}
\left(\begin{matrix}
 \begin{Bmatrix}\displaystyle h_j+\sum\limits_{l=j+1}^{k}\left(\sum\limits_{s=1}^{r_{(k)}}sn_{s}^{l}\bigl(\mathbf{\hat{c}}_{l}^+,-d_{(j)} \mathbf{c}_{j}\bigr)_{D_0R}\right)\vspace{1mm}\\n_0^j,n_{1}^j,\dots,n_{r_{(j)}}^j\end{Bmatrix}
\end{matrix}\prod\limits_{s=1}^{r_{(j)}}z_{i_j,s}^{n_s^j}\right)
 \hat{\mathbf{y}}^{\sum\limits_{j=1}^{k}(\sum\limits_{s=1}^{r_{(j)}}sn_s^j)\mathbf{c}_{j}^+},
 \end{align*}
 where the sum takes over all non-negative integers $n_0^1,n_{1}^{1},\dots,n_{r_{(1)}}^1,\dots,n_0^{k},n_{1}^{k},\dots,n_{r_{(k)}}^k\in\mathbb{Z}_{\ge0}$.

\end{Lemma}
\begin{proof}
 We prove the following claim by induction: for all $1\le p\le k$,
\begin{align}
\prod\limits_{j=p}^{k}L_j^{h_j}
={}&\sum_{n_0^p,n_{1}^{p},\dots,n_{r_{(p)}}^p;\dots;n_0^{k},n_{1}^{k},\dots,n_{r_{(k)}}^k\in\mathbb{Z}_{\ge 0}}\prod_{j=1}^{p-1}L_j^{(\sum\limits_{l=p}^{k}\sum\limits_{s=1}^{r_{(l)}}s n_{s}^{l}(\mathbf{\hat{c}}_{l}^+,-d_{(j)}\mathbf{c}_{j})_{D_0R})}\nonumber\\
&\times \prod\limits_{j=p}^{k}\begin{pmatrix}\!\!
 \begin{Bmatrix}\displaystyle h_j+\sum\limits_{l=j+1}^{k}\left(\sum\limits_{s=1}^{r_{(l)}}s n_{s}^{l}(\mathbf{\hat{c}}_{l}^+,-d_{(j)}\mathbf{c}_{j})_{D_0R}\right)\!\!\vspace{1mm}\\ n_0^j,n_{1}^j,\dots,n_{r_{(j)}}^j\end{Bmatrix}\prod\limits_{s=1}^{r_{(j)}}z_{i_j,s}^{n_s^j}\!
\end{pmatrix}\hat{\mathbf{y}}^{\sum\limits_{j=p}^{k}(\sum\limits_{s=1}^{r_{(j)}}s n_s^j)\mathbf{c}_{j}^+}.\!\!\!\!\label{E38}
 \end{align}
 For any integer $h$ and $1\leq l\leq k$, by Lemma \ref{expansion of polynomial}, we have
 \begin{align*}
 L_l^{h}&=\left(\sum\limits_{s=0}^{r_{(l)}}z_{i_l,s}
 \left(\hat{\mathbf{y}}^{\mathbf{c}_{l}^+}\prod\limits_{j=1}^{l-1}L_j^{(\mathbf{\hat{c}}_{l}^+,-d_{(j)}\mathbf{c}_{j})_{D_0R}}\right)^s\right)^{h} \\ &=\sum\limits_{n_{0},n_1,\dots,n_{r_{(l)}}\in\mathbb{Z}_{\ge 0}}\begin{pmatrix}\begin{Bmatrix}h\\ n_0,n_1,\dots,n_{r_{(l)}}\end{Bmatrix}\prod\limits_{s=1}^{r_{(l)}}\left(z_{i_l,s}\hat{\mathbf{y}}^{s \mathbf{c}_{l}^+}\prod\limits_{j=1}^{l-1}L_j^{s(\mathbf{\hat{c}}_{l}^+,-d_{(j)}\mathbf{c}_{j})_{D_0R}}\right)^{n_s}\end{pmatrix} \\
  &=\sum\limits_{n_{0},n_1,\dots,n_{r_{(l)}}\in \Z_{\geq 0}}\prod\limits_{j=1}^{l-1}L_j^{\sum\limits_{s=1}^{r_{(l)}}s n_s(\mathbf{\hat{c}}_{l}^+,-d_{(j)}\mathbf{c}_{j})_{D_0R}}\begin{Bmatrix}h\\ n_0,n_1,\dots,n_{r_{(l)}}\end{Bmatrix}\left(\prod\limits_{s=1}^{r_{(l)}}z_{i_l,s}^{n_s}\right)\hat{\mathbf{y}}^{\sum\limits_{s=1}^{r_{(l)}}s n_{s}\mathbf{c}_{l}^+}.
 \end{align*}
 By taking $l=k$ and $h=h_k$, this prove \eqref{E38} for $p=k$.

 Now suppose that \eqref{E38} is true for $p+1$ , then multiply both sides by $L_p^{h_p}$, we get
 \begin{align*}
 \prod\limits_{j=p}^{k}L_j^{h_j}
 ={}&L_p^{h_p}\prod\limits_{j=p+1}^{k}L_j^{h_j} \\
={}&\sum_{n_0^{p+1},n_{1}^{p+1},\dots,n_{r_{(p+1)}}^{p+1};\dots;n_0^{k},n_{1}^{k},\dots,n_{r_{(k)}}^k\in \Z_{\geq 0} }\prod_{j=1}^{p-1}L_j^{\bigl(\sum\limits_{l=p+1}^{k}\sum\limits_{s=1}^{r_{(l)}}sn_{s}^{l}(\mathbf{\hat{c}}_{(l)}^+,-d_{(j)}\mathbf{c}_{j})_{D_0R}\bigr)} \\
&\times L_{p}^{h_{p}+\bigl(\sum\limits_{l=p+1}^{k}\sum\limits_{s=1}^{r_{(l)}}sn_{s}^{l}(\mathbf{\hat{c}}_{(l)}^+,-d_{(p)}\mathbf{c}_{p})_{D_0R}\bigr)} \\
&\times \prod\limits_{j=p+1}^{k}\begin{pmatrix}
 \begin{Bmatrix}\displaystyle h_j+\sum\limits_{l=j+1}^{k}\left(\sum\limits_{s=1}^{r_{(l)}}sn_{s}^{l}(\mathbf{\hat{c}}_{l}^+,-d_{(j)}\mathbf{c}_{j})_{D_0R}\right)\vspace{1mm}\\ n_0^j,n_{1}^j,\dots,n_{r_{(j)}}^j\end{Bmatrix}\prod\limits_{s=1}^{r_{(j)}}z_{i_j,s}^{n_s^j}
\end{pmatrix}\hat{\mathbf{y}}^{\sum\limits_{j=p+1}^{k}\bigl(\sum\limits_{s=1}^{r_{(j)}}sn_s^j\bigr)\mathbf{c}_{j}^+}\\
={}&\sum_{n_0^p,n_{1}^{p},\dots,n_{r_{(p)}}^p;\dots;n_0^{k},n_{1}^{k},\dots,n_{r_{(k)}}^k\in \Z_{\geq 0} }\prod_{j=1}^{p-1}L_j^{\bigl(\sum\limits_{l=p}^{k}\sum\limits_{s=1}^{r_{(l)}}s n_{s}^{l}(\mathbf{\hat{c}}_{l}^+,-d_{(j)}\mathbf{c}_{j})_{D_0R}\bigr)} \\
&\times \prod\limits_{j=p}^{k}\begin{pmatrix}\!\!
 \begin{Bmatrix}\displaystyle h_j+\sum\limits_{l=j+1}^{k}\left(\sum\limits_{s=1}^{r_{(l)}}s n_{s}^{l}\bigl(\mathbf{\hat{c}}_{l}^+,-d_{(j)}\mathbf{c}_{j}\bigr)_{D_0R}\right)\!\!\vspace{1mm}\\ n_0^j,n_{1}^j,\dots,n_{r_{(j)}}^j\end{Bmatrix}\prod\limits_{s=1}^{r_{(j)}}z_{i_j,s}^{n_s^j}\!
\end{pmatrix}\hat{\mathbf{y}}^{\sum\limits_{j=p}^{k}\bigl(\sum\limits_{s=1}^{r_{(j)}}s n_s^j\bigr)\mathbf{c}_{j}^+}.\tag*{\qed}
 \end{align*}\renewcommand{\qed}{}
\end{proof}

By applying Lemma \ref{L4} and Theorem \ref{T5}, we have the alternative sum version of Gupta's formula \eqref{E34}.
\begin{Theorem}\label{T6}
 Under the assumption of Theorem {\rm\ref{T5}}, we have
 \begin{align}
 F_{i_k;t_k}(\hat{\mathbf{y}},\mathbf{z})=\sum_{(n_0^1,n_{1}^{1},\dots,n_{r_{(1)}}^1;\dots;n_0^{k},n_{1}^{k},\dots,n_{r_{(k)}}^k)\in\mathbb{Z}_{\ge0}}
 \prod\limits_{j=1}^{k}\left(A_j\prod\limits_{s=1}^{r_{(j)}}z_{i_j,s}^{n_s^j}\right)\hat{\mathbf{y}}^{\sum_{j=1}^{k}(\sum_{s=1}^{r_{(j)}}sn_s^j)\mathbf{c}_{j}^+},\label{E41}
 \end{align}
 where
 \begin{align*}
 A_j=\begin{Bmatrix}\displaystyle -(\mathbf{g}_{k},d_{(j)}\mathbf{c}_{j})_{D_0R}+\sum\limits_{l=j+1}^{k}
 \left(\sum\limits_{s=1}^{r_{(l)}}sn_{s}^{l}(\mathbf{\hat{c}}_{l}^+,-d_{(j)}\mathbf{c}_{j})_{D_0R}\right)\vspace{1mm}\\ n_0^j,n_{1}^j,\dots,n_{r_{j}}^j\end{Bmatrix}.
 \end{align*}
\end{Theorem}
\begin{Remark}
 If $\mathbf{r}=(1,\dots, 1)$, then formulas \eqref{E34} and \eqref{E41} specialize to \cite[Theorems 3.1 and~6.2]{LMN2023}.
\end{Remark}

\subsection*{Acknowledgements}

The authors are grateful to Professor Xueqing Chen for helpful comments. The authors thank the referees for their valuable comments and suggestions in making this article more readable. This work is partially supported by the National Natural Science Foundation of China (Grant No. 11971326).

\pdfbookmark[1]{References}{ref}
\LastPageEnding

\end{document}